\documentclass[reqno,10pt,centertags]{amsart}
\usepackage{amsmath,amsthm,amscd,amssymb,latexsym,esint,upref,stmaryrd,
enumerate,color,verbatim,yfonts}

\usepackage{hyperref} 
\newcommand*{\mailto}[1]{\href{mailto:#1}{\nolinkurl{#1}}}
\newcommand{\arxiv}[1]{\href{http://arxiv.org/abs/#1}{arXiv:#1}}

\newcommand{\Lf}{f_L}

\newcommand{\LSf}{(Sf)_L}
\newcommand{\Lxi}{\xi_L}

\newcommand{\C}{{\mathbb C}}

\newcommand{\bbC}{{\mathbb{C}}}

\newcommand{\bbJ}{{\mathbb{J}}}

\newcommand{\bbN}{{\mathbb{N}}}

\newcommand{\bbR}{{\mathbb{R}}}

\newcommand{\bbZ}{{\mathbb{Z}}}

\newcommand{\bsA}{{\boldsymbol{A}}}

\newcommand{\bsD}{{\boldsymbol{D}}}

\newcommand{\bsH}{{\boldsymbol{H}}}
\newcommand{\bsI}{{\boldsymbol{I}}}

\newcommand{\cB}{{\mathcal B}}

\newcommand{\cF}{{\mathcal F}}

\newcommand{\cH}{{\mathcal H}}
\newcommand{\cI}{{\mathcal I}}

\newcommand{\cL}{{\mathcal L}}
\newcommand{\M}{{\mathcal M}}

\newcommand{\cT}{{\mathcal T}}

\DeclareMathOperator{\rank}{rank}
\DeclareMathOperator{\ran}{ran}
\DeclareMathOperator{\dom}{dom}

\DeclareMathOperator{\tr}{tr}

\DeclareMathOperator*{\sgn}{sgn}

\renewcommand{\Im}{\text{\rm Im}}
\renewcommand{\ln}{\text{\rm ln}}

\newcommand{\loc}{\operatorname{loc}}

\newcommand{\beq}{\begin{equation}}
\newcommand{\enq}{\end{equation}}

\newcommand{\ind}{\operatorname{ind}}
\newcommand{\no}{\notag}
\newcommand{\lb}{\label}
\newcommand{\f}{\frac}

\let\geq\geqslant
\let\leq\leqslant

\makeatletter
\def\theequation{\@arabic\c@equation}

\allowdisplaybreaks
\numberwithin{equation}{section}

\newtheorem{theorem}{Theorem}[section]  
\newtheorem{proposition}[theorem]{Proposition}
\newtheorem{lemma}[theorem]{Lemma}
\newtheorem{corollary}[theorem]{Corollary}
\newtheorem{definition}[theorem]{Definition}
\newtheorem{hypothesis}[theorem]{Hypothesis}
\newtheorem{example}[theorem]{Example}
\theoremstyle{remark}
\newtheorem{remark}[theorem]{Remark}

\begin{document}

\title[The Spectral shift function and the Witten index]{The Spectral shift function and the 
Witten index}

\author[A.\ Carey]{Alan Carey}  
\address{Mathematical Sciences Institute, Australian National University, 
Kingsley St., Canberra, ACT 0200, Australia
and School of Mathematics and Applied Statistics, University of Wollongong, NSW, Australia,  2522}  
\email{\mailto{acarey@maths.anu.edu.au}}
\urladdr{\url{http://maths.anu.edu.au/~acarey/}}
  
\author[F.\ Gesztesy]{Fritz Gesztesy}  
\address{Department of Mathematics,
University of Missouri, Columbia, MO 65211, USA}
\email{\mailto{gesztesyf@missouri.edu}}
\urladdr{\url{http://www.math.missouri.edu/personnel/faculty/gesztesyf.html}}

\author[G.\ Levitina]{Galina Levitina} 
\address{School of Mathematics and Statistics, UNSW, Kensington, NSW 2052,
Australia} 
\email{\mailto{g.levitina@student.unsw.edu.au}}

\author[F.\ Sukochev]{Fedor Sukochev}
\address{School of Mathematics and Statistics, UNSW, Kensington, NSW 2052,
Australia} 
\email{\mailto{f.sukochev@unsw.edu.au}}

\thanks{Submitted to the proceedings of the conference on {\it Spectral Theory and Mathematical Physics, Santiago 2014}, M.\ Mantoiu, G.\ Raikov, and R.\ Tiedra de Aldecoa (eds.), Operator Theory Advances and Application, Birkh\"auser, Springer.}

\date{\today}

\begin{abstract}  
We survey the notion of the spectral shift function of two operators and recent progress on its connection with the Witten index. 
We begin with classical definitions of the spectral shift function $\xi(\, \cdot \,; H_2,H_1)$ under various assumptions on the pair of operators $(H_2, H_1)$ in a fixed Hilbert space and then discuss some of its properties. We then present a  new approach to defining the spectral shift function and discuss Krein's Trace Theorem. In particular, we describe a proof that does not 
use complex analysis \cite{PSZ}  and develop its extension to  general $\sigma$-finite von Neumann algebras $\M$ of type II and unbounded perturbations from the predual of $\M$. 

We also discuss the connection between the theory of the spectral shift function and index theory for certain model operators. We start by introducing various definitions of the Witten index, (an extension of the notion of Fredholm index to non-Fredholm operators). Then we study the model operator $\bsD_\bsA^{} = (d/dt) + \bsA$ in $L^2(\bbR;\cH)$ associated 
with the operator path $\{A(t)\}_{t=-\infty}^{\infty}$, where $(\bsA f)(t) = A(t) f(t)$ for 
a.e.\ $t\in\bbR$, and appropriate $f \in L^2(\bbR;\cH)$ (with $\cH$ being a separable, complex 
Hilbert space). The setup permits the operator family $A(t)$ on $\cH$ to be an unbounded relatively trace class 
perturbation of the unbounded self-adjoint operator $A_-$, and no discrete spectrum 
assumptions are made on the asymptotes $A_{\pm}$. 

When there is a spectral gap for the operators $A_\pm$ at zero, it is shown that the operator $\bsD_\bsA^{}$ is Fredholm and the Fredholm index can be computed as  
\begin{align*}
\ind (\bsD_\bsA^{}) = \xi(0_+; |\bsD_\bsA^{*}|^2, |\bsD_\bsA^{}|^2) = \xi(0; A_+, A_-).
\end{align*}
When $0\in\sigma(A_+)$ (or $0\in\sigma(A_-)$), the operator $\bsD_\bsA^{}$ ceases to be Fredholm. However, under the additional assumption that $0$ is a right and a left Lebesgue 
point of $\xi(\,\cdot\,\, ; A_+, A_-)$, it is proved that $0$ is also a right 
Lebesgue point of $\xi(\,\cdot\,\, ; |\bsD_\bsA^{*}|^2, |\bsD_\bsA^{}|^2)$. For the resolvent (resp., semigroup) regularized Witten index $W_r(\bsD_\bsA^{})$ (resp., $ W_s(\bsD_\bsA^{})$) the following equality holds,  
\begin{align*} 
W_r(\bsD_\bsA^{}) &= W_s(\bsD_\bsA^{})= \xi(0_+; |\bsD_\bsA^{*}|^2, |\bsD_\bsA^{}|^2)       \\ 
& = [\xi(0_+; A_+,A_-) + \xi(0_-; A_+, A_-)]/2.
\end{align*}

We also study a special example, when the perturbation of the unbounded self-adjoint operator $A_-$ is not assumed to be relatively trace class. In this example $A_-= - i \frac{d}{dx}$ is the differentiation operator on $L^2(\bbR)$ and the perturbation is given by the multiplication operator by a (bounded) real-valued function $f$ on $\bbR$. Under certain assumptions on $f$ it is proved   that 
\begin{align*} 
W_r(\bsD_\bsA^{}) &= W_s(\bsD_\bsA^{})=  \xi(0_+; \bsD_\bsA^{}\bsD_\bsA^{*}, \bsD_\bsA^{*}\bsD_\bsA^{})\\&= \xi(0; A_+,A_-)=\frac{1}{2\pi}\int_\bbR f(s) \, ds.
\end{align*} 
\end{abstract}

\maketitle

\newpage 

{\scriptsize{\tableofcontents}}

\section{Introduction}

The purpose of this article is twofold:  We give a detailed survey of the Lifshitz--Krein spectral shift function and its properties, and we 
 then review the notion of the Witten index and its relation to the spectral shift function and to spectral flow. 

We begin in section 2 with  an account of the history of the spectral shift function starting with the work of Lifshitz and Krein. We discuss several points of view
on the definition and then move on to more recent developments. 
We explain in some detail a recent real analysis approach to the
fundamental  theorem of Krein (almost all complete earlier proofs use complex analysis, see, however, \cite{Si94} and \cite{Vo87}).
The novelty here is that the proof also applies when one works in the generality of semifinite von Neumann algebras (rather than just the 
algebra of bounded operators on a Hilbert space).

Starting in section 3 we survey the properties of the Witten index from a more contemporary perspective. We introduce a special ``supersymmetric'' model operator motivated by geometric considerations.  We describe in section 4 recent results relevant to index theory that do not depend on assuming that the operators under study all have discrete spectrum.  In particular, we focus on two formulae (we call these the principle
trace formula and the Pushnitski formula) that seem especially interesting.
Generalisations of both of these formulae are described in terms of recent results (published and, as of yet, unpublished ones).  We briefly explain in the final section some new
examples that point the way to higher-dimensional examples.

\section{Spectral shift function}

In 1947, the well-known physicist Lifshitz considered perturbations of an operator
 $H_0$ (arising as the Hamiltonian of a lattice model in quantum mechanics)  by a finite-rank  perturbation $V$ and found some formulae
and quantitative relations for the size of the shift of the eigenvalues.
In one of his papers the spectral shift function (SSF), $\xi(\, \cdot \,; H_0+V, H_0)$, appeared
for the first time, and formulae for it in the case of a finite-rank perturbation
were obtained. 

Lifshitz later continued these investigations and applied them
to the problem of computing the trace of the operator $\phi(H_0+V)-\phi(H_0)$, where $H_0$ is the  unperturbed self-adjoint operator, $V$ is a self-adjoint, finite-dimensional perturbation,
and $\phi$ is an appropriate function (belonging to a fairly broad class). He obtained (or, rather, surmised) the remarkable relation 
 \begin{equation}\label{Krein_tr_formula}
 \tr(\phi(H_0+V)-\phi(H_0))=\int_\bbR \phi'(\lambda)\xi(\lambda; H_0+V, H_0) \, d\lambda,
 \end{equation}
where the function $\xi(\, \cdot \,; H_0+V, H_0)$ depends on operators $H_0$ and $V$ only.

A physical example treated by Lifshitz is the following: if $H_0$ is the operator describing the oscillations
of a crystal lattice, then the free energy of the oscillations can be represented in
the form $F=\tr(\phi(H_0))$,
for some $\phi$. In this case, the trace formula enables one to compute the change in the free energy of oscillations of the crystal
lattice upon introduction of a foreign admixture into the crystal.

If one wants to study continuous analogues of lattice models, perturbations $V$, as a rule, are no longer described by finite-rank operators. For such models the appropriate class of perturbations, such that the spectral shift function may be defined, needs to be described. 
In his paper \cite{Kr53}, M.~G.~Krein resolved this problem. Furthermore, he described the broad class of functions $\phi$ for which \eqref{Krein_tr_formula} holds.
His approach was based on the notion of perturbation determinants to be discussed next.

\subsection{Perturbation determinants.}
Let $\cH$ be a complex, separable Hilbert
space, $\mathcal{L}(\cH)$ be the algebra of all
bounded linear operators in $\cH$ and let $\cL_1(\cH)$ be the ideal of all trace class operators.
The latter ideal, besides carrying the standard trace, also gives rise to the notion of a determinant, which generalizes the corresponding notion in the finite-dimensional case. Let $T\in\cL_1(\cH)$. For any orthonormal basis $\{\omega_n\}_{n\in\bbN}$ in $\cH$ consider the $N\times N$ matrix $\cT_N$ with elements $\delta_{m,n}+(T\omega_m,\omega_n)$, $m,n \in 1,\dots, N$. Then the following limit exists 
$$\lim_{N\rightarrow\infty}\det(I+\cT_N)=:\det(I+T),$$
independently of the choice of the basis $\{\omega_n\}_{n\in\bbN}$ (cf., \cite[Ch.~IV]{GK69}).
The functional $\det (I + \, \cdot \,):\cL_1(\cH)\to \bbC$ is called the \emph{determinant}; it is continuous with respect to  the $\cL_1(\cH)$-norm.

In terms of eigenvalues of $T\in \cL_1(\cH)$, $\{\lambda_k(T)\}_{k \in \cI}$, $\cI \subseteq \bbN$, an appropriate index set, one has 
$$\det(I+T)=\prod_{k \in \cI} (1+\lambda_k(T)),$$
where the product converges absolutely (due to the fact that $\sum _{k \in \cI} |\lambda_k|<\infty$).
We note the following properties of the determinant \cite{GK69}
\begin{gather*}
\det(I+T^*)=\overline{\det(I+T)},\quad T\in\cL_1(\cH)\\
\det(I+T_1)(I+T_2)=\det(I+T_1)\det(I+T_2),\quad T_1,T_2\in\cL_1(\cH)\\
\det(I+T_1T_2)=\det(I+T_2T_1),\quad T_1, T_2 \in \cL(\cH), \; T_1T_2, T_2T_1\in\cL_1(\cH).
\end{gather*}

In the following, let $H_0$, $H$ be self-adjoint operators in $\cH$ with $\dom(H_0)=\dom(H)
$, and let $V=H-H_0$. Assume that $VR_z(H_0)\in\cB_1(\cH),$
where $R_z(T)$ denotes the resolvent of an operator $T$, that is, $R_z(T)=(T-z I)^{-1}$.
Under these assumptions one can introduce the \emph{perturbation determinant}
$$\Delta(z)=
\Delta_{H/H_0}(z):=\det(I+VR_z(H_0)) = \det\big((H - z I)(H_0 - z I)^{-1}\big), \quad\Im(z)\neq 0.$$
Next we briefly recall some properties of perturbation determinants.

For self-adjoint operators $H_0, H$ the mapping $z\to\Delta_{H/H_0}(z)$ is analytic in both the half-planes $\Im(z)>0$ and $\Im(z)<0$ and
$$\Delta_{H/H_0}(\bar{z})=\overline{\Delta_{H/H_0}(z)}, \quad \Im(z)\neq 0.$$
One has $\Delta_{H/H_0}(z)\neq 0$ for $ \Im(z)\neq 0$.

In addition, since $V\in \cB_1(\cH)$, standard properties of resolvents imply that
$$\|VR_{H_0}(z)\|_1\rightarrow 0 \, \text{ as } \, |\Im(z)|\rightarrow\infty,$$
and therefore, 
$$ \Delta_{H/H_0}(z)\rightarrow 1 \, \text{ as $ |\Im(z)|\rightarrow\infty$.}$$
Since the function $\Delta_{H/H_0}(\cdot)$ is analytic in the open upper and lower half plane and 
since $\Delta_{H/H_0}(z)\neq 0, \, \Im(z)\neq 0$, it is a standard fact from complex analysis that there exists a function $G(\cdot)$ analytic in both of the upper and lower half planes such that
$e^G=\Delta_{H/H_0}.$
Naturally, one denotes the function $G$ by $\ln(\Delta_{H/H_0})$. It is clear that the function 
$\ln(\Delta_{H/H_0})$ is multivalued and its different values at a point $z$, $\Im(z)\neq 0$, differ by $2\pi i k$, $k \in \bbZ$.
Since $\Delta_{H/H_0}(z)\rightarrow 1, $ as $|\Im(z)|\rightarrow\infty,$ one \emph{fixes the branch} of the function $\ln(\Delta_{H/H_0})$ by requiring that
$\ln(\Delta_{H/H_0}(z)) \rightarrow 0$ as $ |\Im(z)|\rightarrow\infty.$

\subsection{Construction of the SSF due to M.~G.~Krein.}
To construct the spectral shift function by Krein's method we exploit the following representation of the function 
$\ln(\Delta_{H/H_0}(z))$, 
\begin{equation}\label{repr_lnDelta}
\ln(\Delta_{H/H_0}(z))=\int_\bbR \f{\xi(\lambda;H,H_0) \, d\lambda}{\lambda-z}, \quad \Im(z)\neq 0,
\end{equation}
with a real-valued $\xi(\, \cdot \,; H,H_0) \in L_1(\bbR)$.

The proof of \eqref{repr_lnDelta} relies on the following classical result from complex analysis.

\begin{theorem}[Privalov representation theorem]\label{Privalov_thm}
Suppose that $F$
is holomorphic in the open upper half-plane. If $\Im(F)$ is bounded and non-negative $($respectively, non-positive\,$)$ and if $\sup_{y \geq 1} y|F(iy)|<\infty,$
then there exists a nonnegative $($respectively, non-positive\,$)$ real-valued function $\xi\in L_1(\bbR)$
such that
$$F(z)=\int_\bbR\frac{\xi(\lambda) \, d\lambda}{z-\lambda}, \quad \Im(z)>0.$$ The function $\xi$ is uniquely determined by the Stieltjes inversion formula, 
$$\xi(\lambda)=\frac1\pi\lim_{\varepsilon\downarrow0+}\Im (F(\lambda+i\varepsilon)) \, 
\text{ for a.e.~$\lambda\in\bbR.$}$$
\end{theorem}

Next we sketch the proof of the first theorem of Krein (see Theorem \ref{thm_Krein}).

To verify the assumptions in Privalov's Theorem for $F=\ln(\Delta_{H/H_0})$,  Krein proceeded as follows:
\begin{itemize}
\item First, suppose that $\rank(V)=1,$ that is, $V=\gamma(\cdot, h)h,$ \, $h\in\cH,$ \, $\|h\|=1,$ \, $\gamma\in\bbR$.  Then
$$\Delta_{H/H_0}(z)=1+\gamma(R_{H_0}(z)h,h).$$

\end{itemize}

Using this explicit form of the perturbation determinant one can prove that the function 
$\ln(\Delta_{H/H_0}(\cdot))$ satisfies all the assumptions in Privalov's theorem (for details see, e.g., Yafaev's book \cite{Ya92}).
Hence, there exists a function $\xi(\lambda; H,H_0)$ satisfying \eqref{repr_lnDelta}, and furthermore, the function $\xi(\, \cdot \,; H,H_0)$ can be expressed in the form
\begin{equation}\label{repr_xi}
\xi(\lambda; H,H_0)=\frac1\pi\lim_{\varepsilon\rightarrow+0} \Im(\ln(\Delta_{H/H_0}(\lambda+i\varepsilon))),\quad \text{a.e. } \lambda\in\bbR.
\end{equation}

\begin{itemize}
\item Suppose now, that $\rank(V)=n<\infty,$ that is, 
$$V=\sum_{k=1}^n\gamma_k(\cdot,h_k)h_k, \quad \gamma_k=\bar{\gamma}_k,\, \|h_k\|_1, \, 1\leq k\leq n.$$
Denoting $$V_m=\sum_{k=1}^m\gamma_k(\cdot,h_k)h_k,\quad H_m=H_0+V_m, \quad 1\leq m\leq \rank(V),$$
one infers that the difference $H_m-H_{m-1}$ is a rank-one operator. In addition, by the multiplicative property of the determinant one concludes that
\begin{equation}\label{Delta_n}
\ln(\Delta_{H/H_0}(z)) = \sum_{m=1}^n\ln(\Delta_{H_m/H_{m-1}}(z)).
\end{equation}
Applying the first step to the operators $H_m, H_{m-1}$ one infers the existence of the corresponding SSFs $\xi(\, \cdot \,; H_m,H_{m-1})$, $1 \leq m \leq \rank(V)$.

\end{itemize}

Set
$$\xi(\lambda;H,H_0)=\sum_{m=1}^n\xi(\lambda; H_m,H_{m-1}), \quad 1\leq k\leq n.$$ 
There are $L_1(\bbR)$-norm estimates for each $\xi(\, \cdot \,; H_m,H_{m-1})$  which  ensure that the function $\xi(\lambda;H,H_0)$ is integrable. Furthermore, since for every $m$, the representations \eqref{repr_lnDelta} and \eqref{repr_xi} for  $\ln(\Delta_{H_m/H_{m-1}})$ and 
$\xi(\lambda; H_m,H_{m-1})$, respectively, hold, one can infer from \eqref{Delta_n} and the definition of $\xi(\lambda;H,H_0)$ that representations \eqref{repr_lnDelta} and \eqref{repr_xi} hold also for $\ln(\Delta_{H/H_0})$ and $\xi(\lambda;H,H_0)$.

\begin{itemize}
\item Suppose now, that $V$ is an arbitrary trace class perturbation . Let $V_n$ be a sequence of finite-rank operators, such that $\|V-V_n\|_1\rightarrow0, \, n\rightarrow \infty$. Set
$$\xi(\lambda;H,H_0)=\sum_n \xi(\lambda; H_n,H_{n-1}),$$
where the sum now is infinite (unless, $V$ is a finite-rank operator).

Then, convergence properties of determinants and the $L_1(\bbR)$-norm estimate for each 
$\xi(\, \cdot \,; H_n,H_{n-1})$ imply that this series converges in $L_1(\bbR)$ and all the desired representations \eqref{repr_lnDelta} and \eqref{repr_xi} for  $\ln(\Delta_{H/H_0})$ and 
$\xi(\, \cdot \,;H,H_0)$ hold.
\end{itemize}

The following result is the first theorem of M.~G.~Krein.

\begin{theorem}\cite{Kr53}\label{thm_Krein}
Let $V\in\cB_1(\cH)$ be self-adjoint. Then the following representation holds
$$\ln(\Delta_{H/H_0}(z)) = \int_\bbR \f{\xi(\lambda; H,H_0) \, d\lambda}{\lambda-z}, \quad \Im(z)\neq 0,$$
where
\begin{equation}\label{Krein}
\xi(\lambda; H,H_0)=\frac1\pi\lim_{\varepsilon\downarrow 0}\Im(\ln(\Delta_{H/H_0}(\lambda+i\varepsilon))) 
\, \text{ for a.e.~$\lambda \in \bbR$,}
\end{equation}
in particular, the limit in \eqref{Krein} exists for a.e. $\lambda\in\bbR$. In addition, 
\begin{equation}\label{estim_xi_trclass}
\int_\bbR|\xi(\lambda; H,H_0)| \, d\lambda\leq \|V\|_1,\quad 
\int_\bbR\xi(\lambda; H,H_0) \, d\lambda = \tr(V).
\end{equation}
Moreover, $\xi(\lambda; H,H_0)\leq k_+$ $($respectively, $\xi(\lambda;H,H_0)\geq -k_-$$)$ for a.e. $\lambda\in \bbR$, provided that the perturbation $V$ has only $k_+$ positive $($respectively, $k_-$ negative\,$)$ eigenvalues.
\end{theorem}

Next, we turn to the rigorously proved trace formula, which is now customarily  referred as the Lifshitz--Krein trace formula.

\begin{theorem}[Second theorem of M.~G.~Krein]\label{Krein_2}
Let $V\in\cB_1(\cH)$ and assume that $f\in C^1(\bbR)$ and its derivative admits the representation
$$f'(\lambda)=\int_\bbR\exp(-i\lambda t) \, dm(t), \quad |m|(\bbR)<\infty,$$
for a  finite $($complex\,$)$ measure $m$. Then
$[f(H)-f(H_0)] \in\cB_1(\cH),$ and the following trace formula holds
\begin{equation}\label{trace_formula}
\tr(f(H)-f(H_0))=\int_\bbR f'(\lambda) \xi(\lambda;H,H_0) \, d\lambda.
\end{equation} 
\end{theorem}

\begin{remark}
$(i)$ As is clearly seen from the arguments sketched, Krein's original proof was based on complex analysis. Attempts to produce a ``real-analytic proof'' are discussed later. \\[1mm] 
$(ii)$ The function $\xi(\, \cdot \,; H, H_0)$ is an element of $L_1(\bbR)$, that is, it represents an \emph{ equivalence class} of Lebesgue measurable functions. Therefore, generally speaking, the notation $\xi(\lambda; H, H_0)$ is meaningless for a fixed $\lambda\in\bbR$. \\[1mm] 
$(iii)$ For a trace class perturbation $V$, the spectral shift function $\xi(\, \cdot \,; H,H_0)$ is \emph{unique}.  \\[1mm] 
$(iv)$ The Lifshitz--Krein trace formula can be extended in various ways. One could attempt to describe the class of functions $f$, for which this formula holds; however, we will not cover this direction. Another important direction is to enlarge the class of perturbations $H-H_0$. We shall present some results in this direction below. \hfill $\diamond$
\end{remark}

\subsection{Properties of the spectral shift function.} \label{sec_properties_ssf}
Let $H_0,H_1$ and $H$ be such that $(H_1-H_0), (H-H_1) \in\cB_1(\cH)$.
First, we will list the simplest properties of the SSF. 

These are that for a.e. $\lambda\in\bbR$ we have
$$\xi(\lambda;H,H_1)+\xi(\lambda; H_1,H_0)=\xi(\lambda; H,H_0),$$
in particular, $\xi(\lambda; H,H_0)=-\xi(\lambda; H_0,H)$, and also
the inequality
$$\|\xi(\, \cdot \,; H,H_0)-\xi(\, \cdot \,; H_1,H_0)\|_1\leq\|H-H_1\|_1$$ holds. 
 In addition, if $H\geq H_1$, then
$$\xi(\lambda;H,H_0)\geq\xi(\lambda; H_1,H_0) \, \text{ for a.e. $\lambda\in\bbR$.}$$

Next, we describe some special situations where one can select concrete representatives from the equivalence class $\xi(\, \cdot \,; H, H_0)$, which justifies the term``the spectral shift function". These properties of the SSF $\xi(\, \cdot \,; H,H_0)$ are  associated with the spectra of the operators $H_0$ and $H$. For the complete proof we refer to \cite[Ch.~8]{Ya92}

\begin{enumerate}

\item[$(i)$] Let $\delta$ be an interval (possibly unbounded) such that $\delta\subset \rho(H_0)\cap\rho(H)$. Then $\xi(\, \cdot \,; H,H_0)$ takes a \emph{constant integer value} on $\delta$, that is, 
$$\xi(\lambda;H,H_0)=n,\quad n\in\bbZ,\;  \lambda \in\delta.$$
If the interval $\delta$ contains a half-line, then the $L^1$-condition on $\xi$ implies that $n=0$.

\item[$(ii)$] Let $\mu$ be an isolated eigenvalue of multiplicity $\alpha_0<\infty$ of $H_0$ and multiplicity $\alpha$ for $H$.  Then
\begin{equation}\label{xi_one_eigenvalue}
\xi(\mu_+; H,H_0) - \xi(\mu_-; H,H_0)=\alpha_0-\alpha.
\end{equation}

\end{enumerate}

 Property $(ii)$ can be generalized as follows: 
\begin{enumerate}
\item[$(iii)$] Suppose that in some interval $(a_0,b_0)$ the spectrum of $H_0$ is discrete (i.e., the spectrum of $H_0$ consists at most of eigenvalues of $H_0$ of \emph{ finite } multiplicity all of which are \emph{isolated} points of $\sigma(H_0)$). Then, by Weyl's theorem on the invariance of essential spectra (see, e.g., \cite[Theorem 5.35]{Ka80}), $H$ has discrete spectrum in $(a_0,b_0)$ as well.

Let $\delta=(a,b),\, a_0<a<b<b_0.$
Introduce the \emph{eigenvalue counting functions} $N_0(\delta)$ and  $N(\delta)$ of the operators $H_0$ and $H$, respectively, in the interval $\delta$ as the sum of the multiplicities of the eigenvalues in $\delta$ of the operator $H_0$, respectively, $H$. Since the interval $\delta$ is finite and both operators $H_0,H$ have discrete spectrum, $N_0(\delta)$ and  $N(\delta)$ are finite.
In this case one has the equality, 
\begin{equation}\label{xi_on_discrete}
\xi(b_-; H,H_0)-\xi(a_+; H,H_0)=N_0(\delta)-N(\delta).
\end{equation}
\end{enumerate}

The preceding property implies, in particular, the following fact. 
\begin{enumerate}
\item[$(iv)$] Let $H_0$ be a nonnegative self-adjoint operator with purely discrete spectrum (i.e., 
$\sigma_{ess}(H_0) = \emptyset$). Since the perturbation $V$ is trace class, there exists $c\in\bbR$, such that $H\geq c,$ that is, $H$ is also lower semibounded. Generally, $H$ will of course not be nonnegative and so one should expect negative eigenvalues of $H$. Thus, property $(iii)$ implies that for $\lambda<0$, 
$$\xi(\lambda_-)=-N(\lambda,H),$$
where $N(\lambda,H)$ is the sum of multiplicities of the eigenvalues of $H$ lying to the left of the point $\lambda<0$.
\end{enumerate}

On the other hand, the following result demonstrates that any function from $L_1(\bbR)$ arises as the spectral shift function for some pair of operators.
\begin{enumerate}
\item[$(v)$] Let $\xi$ be an arbitrary real-valued element of $L_1(\bbR)$. Then, there exists a pair of self-adjoint operators $H_0,H$, such that $(H-H_0) \in\cB_1(\cH)$ and $\xi$ is the SSF $\xi(\, \cdot \,; H,H_0)$ for the pair $(H,H_0)$. In addition, if $0\leq\xi\leq 1$, then there is a pair $H_0,H$ such that $H-H_0$ is a positive rank-one operator \cite{Kr53}, \cite{KY81}.

\end{enumerate}

\subsection{Earlier real-analytic approaches.}
In the following we discuss other approaches for constructing the SSF.
The first attempt to prove the existence of the SSF
without relying on complex analysis was made by
Birman and Solomyak in \cite{BS75}.
This method is based on consideration of the \emph{family of operators,}
$$H_s=H_0+sV, \quad s\in[0,1], \quad H =H_1,$$ 
and their family of spectral measures
$\{E_{H_s}(\lambda)\}_{\lambda \in \bbR}$. 
Employing the \emph{theory of double operator
integrals} also developed by these authors, it can be proved that for sufficiently large class of functions $f$, there exists a continuous derivative in $\cB_1(\cH)$-norm of the operator-valued function $s\mapsto f(H_s)$, represented in the double operator integral form as
$$\frac{df(H_s)}{ds}=\int_\bbR\int_\bbR\frac{f(\mu)-f(\lambda)}{\mu-\lambda} \, dE_{H_s}(\mu)VdE_{H_s}(\lambda).$$

Furthermore, Birman and Solomyak obtained the equality
$$\tr\bigg(\frac{df(H_s)}{ds}\bigg)=\int_\bbR f'(\lambda) \, d\tr(V E_{H_s}(\lambda)).$$
Integration with respect to $s$ then yields the formula
$$\tr(f(H)-f(H_0))=\int_\bbR f'(\lambda) \, d\Xi_{H,H_0}(\lambda),$$
where the \emph{spectral averaging measure $\Xi_{H,H_0}$} is defined by
$$\Xi_{H,H_0}(X)=\int_0^1\tr(V E_{H_s}(X))\,ds,$$
with $X \subseteq \bbR$ a Borel set. 

However, this attempt to yield an alternative proof of Krein's Theorem \ref{Krein_2} was unsuccessful since the authors failed to establish the absolute
continuity of the latter measure with respect to  Lebesgue measure.

We note, that \emph{if} one introduces $\xi(\cdot; H,H_0)$ by Krein's Theorem \ref{thm_Krein}, then
$$\int_X \xi(\lambda;H,H_0) \, d\lambda=\Xi_{H,H_0}(X),$$
for any Borel set $X \subseteq \bbR$, that is, the measure $\Xi$ is indeed absolutely continuous.

The second attempt to deliver a real-analytic proof was due to  Voiculescu \cite{Vo87}, 
his method was based on the classical Weyl--Berg--von Neumann
theorem. However, his attempt also failed to recover the full generality of
Krein's original result.

Another attempt to obtain a proof of
Krein's formula without appealing to complex-analytic methods was introduced by Sinha and Mohapatra \cite{Si94}. Again, that attempt did not yield the full generality of the result  and
does not seem to apply to general semifinite von Neumann algebras.

\subsection{The case of semifinite von Neumann algebras.}
Some problems in noncommutative geometry require replacing the algebra $\cB(\cH)$ of all bounded linear operators on a Hilbert space $\cH$ and unbounded operators on $\cH$ with a general semifinite von Neumann algebra $\M$ and unbounded operators affiliated with $\M$.
A typical example of differential operators affiliated to semifinite von Neumann algebras arises in the context of Atiyah's $L^2$-index theorem and its extensions. (For example, the paper 
\cite{BC15} considers the case of lifts of Dirac-type operators acting on sections of a finite dimensional vector bundle over a complete Riemannian manifold $M$ to a Galois cover 
$\widetilde M$  of $M$.)

The first attempt to extend Krein's results and methods to the realm of semifinite von Neumann algebras was made in \cite{ADS06}. It broadly followed Krein's complex analysis proof. However, it does not offer an adequate extension to general semifinite von Neumann algebras of the notion of the perturbation determinant, which plays the key role in Krein's proof. This difficulty is circumvented in \cite{ADS06} via the use of the notion of a \emph{Brown measure} \cite{Br86}.

The core of the approach in \cite{ADS06} is to show that there exists a neighbourhood of the spectrum of the operator $R_z(H_0)V$, which does not  intersect the half-line $(-\infty,-1]$, in the case where $V\geq 0$ or $-V \geq 0$. One then applies one of the principal results of Brown \cite{Br86} to establish estimates needed for the application of the Privalov representation theorem (see Theorem \ref{Privalov_thm}). Finally, the proof in \cite{ADS06} proceeded under the additional assumption that $H-H_0$ is a bounded trace class perturbation.

Another subsequent paper \cite{ACS07} employed the double operator  integral (DOI) technique due to  Birman and Solomyak, but in a slightly different form suitable for semifinite von Neumann algebras using an approach from \cite{ACDS09}. Following the idea of Birman and Solomyak, one can define the spectral shift measure for a pair $(H,H_0)$, by setting
$$\Xi_{H,H_0}(X)=\int_0^1\tau(VE_{H_s}(X))ds,$$ 
where $\tau$ is a faithful normal semifinite trace on $\M$.
Assuming that $H_0$ has \emph{$\tau$-compact resolvent}, and the perturbation $V$ is \emph{bounded}, it can be proved that the spectral shift measure $\Xi_{H,H_0}$ is \emph{absolutely continuous} with respect to the Lebesgue measure and the resulting Radon--Nikodym derivative is the SSF for the pair $(H,H_0)$.

The first complete ``real analytic proof" of the Lifshitz--Krein formula is due to Potapov, Sukochev,   and Zanin \cite{PSZ}.
That paper delivers a  rather short and straightforward proof of the Lifshitz--Krein formula  without any use of complex analytic tools. The approach in  \cite{PSZ} can be characterized as a combination of methods drawn
from the double operator integration theory of Birman and Solomyak and from Voiculescu's ideas based on the Weyl--Berg--von Neumann theorem. 
The result holds for an arbitrary semifinite von Neumann algebra $\M$, equipped with a faithful normal semifinite trace $\tau$ and (unbounded) operators $H_0,H$ affiliated with $\M$, such that $H-H_0$ belongs to the space $\cL_1(\M,\tau)$, the predual of the algebra $\M$. 

We denote by \emph{$W_1$} the class of all differentiable
functions $f:\mathbb R\to \mathbb R$ such that $\mathcal{F}(f')\in
L_1(\bbR),$ where the symbol $\cF$ denotes the standard Fourier transform. The following theorem is the main result of \cite{PSZ}: 

\begin{theorem}\label{thm_PSZ}
  Let~$\mathcal M$ be a von Neumann algebra equipped with a
  faithful normal semifinite trace~$\tau$.  If the self-adjoint operators
$H_0,H$ affiliated with~$\mathcal{M}$ are such that~$(H-H_0) \in 
\cL_1(\M,\tau)$, then there is a
  function~$\xi(\, \cdot \,; H,H_0) \in L_1(\bbR)$ such that the trace formula
  \begin{equation}
     \tau (f(H)- f(H_0)) =
    \int_{\mathbb R} f'(\lambda)\, \xi(\lambda; H,H_0) \, d\lambda.  \end{equation}
    holds for all $f \in W_1$.
\end{theorem}

\begin{remark}
If the von Neumann algebra $\M$ is the type I factor $\cB(\cH)$ with the standard trace, then Theorem \ref{thm_PSZ} delivers an alternative proof of Krein's result $($i.e., 
Theorem \ref{Krein_2}$)$. \hfill $\diamond$ 
\end{remark}

Below we outline  the proof of Theorem \ref{thm_PSZ}. 

We start by introducing the \emph{distribution function} $N_{H_0}$ of the operator $H_0$, that is, 
$$N_{H_0}(t):=\tau(E_{H_0}(t,\infty)), \quad t\geq 0,$$ where
$E_{H_0}((t,\infty))$ is the spectral projection of the
self-adjoint operator $H_0$ corresponding to the
interval~$(t,\infty).$

The proof in \cite{PSZ} is divided into several stages. For simplicity we denote by 
$\xi^{(j)}(\, \cdot \,;H,H_0)$ the function constructed on the $j$-th step.

{\bf Step~$\boldsymbol{(i)}$.}
 Let the trace $\tau$ be \emph{finite}, that is, $\tau(I)<\infty$ and \emph{$H_0,H\in \mathcal M$}.
In this case, the SSF is merely defined as
\emph{$$\xi^{(1)}(\, \cdot \,;H,H_0) = N_{H}(\cdot) - N_{H_0}(\cdot).$$}
Since the trace $\tau$ is finite, both $N_{H}$ and $N_{H_0}$ are finite.

One should note the similarity of this formula with property $(iii)$ of the SSF (see \eqref{xi_on_discrete}). One can think of this equation as the ``naive'' definition of the SSF. However, while this definition is correct for finite von Neumann algebras, there are examples of self-adjoint operators $H,H_0$ in infinite dimensional Hilbert space with $H-H_0$ being a rank-one operator such that  the operator
$E_{H_0}((t,\infty))  -E_{H}((t,\infty))$
\emph{is not} a trace class operator for all $t$ on the spectrum \cite{Kr53} (the example concerns 
self-adjoint resolvents of Dirichlet and Neumann Laplacians on a half-line).

 The
    function~$\xi^{(1)}(\, \cdot \,;H,H_0)$ is supported
    on the interval $[-a,a]$, where $a:=\max\{\|H_0\|_{\infty},\|H\|_{\infty}\}$.
 Furthermore, it possesses  a property similar to that of the Krein' SSF (see \eqref{estim_xi_trclass}), 
 \begin{equation} 
 \|N_{H}-N_{H_0}\|_{\infty}\leq\tau({\rm supp}(H-H_0)), \quad 
  \|N_{H}-N_{H_0}\|_1\leq\|H-H_0\|_1.      \label{maj lemma1}
  \end{equation}

{\bf Step~$\boldsymbol{(ii)}$.}
In the second step, the trace formula is proved for \emph{bounded} operators $H_0,H\in \mathcal M,$ with the perturbation $V=H-H_0$ being a \emph{nonnegative} operator with \emph{$\tau$-finite support}, and for functions of the form \emph{$f(s)=s^m$}. Here we use an  idea  noted by Voiculescu, who proved the Krein trace formula for the case of polynomials.

Proving a result similar to the classical Berg--Weyl--von Neumann theorem we construct  a
family of~$\tau$-finite projections~$p_n$, $n\in \bbN$, with
$p_n\uparrow I$ such that
\begin{equation}\label{some_eq}
\tau((p_nHp_n)^m-(p_nH_0p_n)^m) - \tau(H^m -H_0^m) \to 0 \text{ as }
n\to \infty.
\end{equation}
Since for every~$n\in\bbN$, $\tau(p_n 1 p_n)<\infty$, by Step~$(i)$, there exists a
positive function~$\xi_n^{(1)}=\xi^{(1)}(\, \cdot \,; p_nHp_n,p_nH_0p_n)$, 
supported on~$[-a,a]$, satisfying the trace formula.
In addition, by \eqref{maj lemma1}, the sequence $\big\{\xi_n^{(1)}\big\}_{n\in \bbN}$ is bounded in $L_\infty((-a,a)).$ By the Banach--Alaoglu Theorem the latter is compact in the weak$^*$-topology, and therefore, there exists a directed set $\bbJ$ and a mapping $\psi:\bbJ \to \bbN$ such that for every $n\in\bbN$, there exists $j(n)\in \bbJ$ such that $\psi(j)>n$ for $j>j(n)$ and such that the 
net $\xi_{\psi(j)}^{(1)}(\, \cdot \,; p_{\psi(j)} H p_{\psi(j)}, p_{\psi(j)} H_0 p_{\psi(j)})$ converges in weak$^*$-topology. The function $\xi^{(2)}(\, \cdot \,\; H, H_0)$ is then defined by 
$${{\xi^{(2)}(\, \cdot \,; H,H_0)}}:=\lim_{j\in\mathbb{I}} 
\xi_{\psi(j)}^{(1)}(\, \cdot \,; p_{\psi(j)} H p_{\psi(j)}, p_{\psi(j)} H_0 p_{\psi(j)}), $$ 
and proved to be the SSF.

{\bf Step~$\boldsymbol{(iii)}$.} Let $H,H_0\in \mathcal M,$ $f\in C^2_b(\mathbb R)$. 
In this step we remove the assumptions $H\geq H_0$   and $\tau({\rm
supp}(H-H_0))<\infty$.

We prove that, without loss of generality, one can assume that $H\geq H_0$. Let~$0\leq D_n\leq
H-H_0$, $n \in \bbN$, be such that~$D_n\uparrow H-H_0$ as $n \to \infty$, and
  $\tau({\rm supp}(D_n))<\infty$, $n\in\bbN$.

Since polynomials are dense in~$C^2([-a,a])$, it follows from Step~$(ii)$ and DOI techniques that 
\begin{equation}\label{tbd2}
\tau(f(H_0+D_n)-f(H_0))=\int_{-a}^a f'(\lambda)\xi^{(2)}(\lambda; H_0+D_n,H_0) \, d\lambda, 
\quad f\in C^2_b(\mathbb{R}). 
\end{equation} 
Then, proving that the sequence  $\big\{\xi^{(2)}(\, \cdot \,; H_0+D_n,H_0)\big\}_{n\in\bbN}$ increases and is uniformly bounded, one infers from the Monotone Convergence
Principle that the sequence $\{\xi^{(2)}(\, \cdot \,; H_0+D_n,H_0)\}_{n\in\bbN}$
converges in $L^1(\bbR)$; its limit is denoted by $\xi^{(3)}(\, \cdot \,; H,H_0).$ This function is now the SSF for the pair $(H,H_0)$.

{\bf Step~$\boldsymbol{(iv)}$.}
  The final step in this approach consists in removing the assumption that the
operators~$H_0$ and~$H$ are bounded. This is the key point of the
proof in which DOI techniques are used in its full strength. This part of the proof is rather technical. We briefly outline the main ideas. 

Choose a $C^2$-bijection  $h: \bbR \to (a,b)$ for some $a, b \in\bbR$, $a < b$. Then 
  the operators $h(H_0)$ and $h(H)$ are bounded, so that  applying
Step~$(iii)$ to the operators~$h(H_0)$ and~$h(H)$, one  defines
  $$\xi^{(4)}(\, \cdot \,; H,H_0):=\xi^{(3)}(\, \cdot \,; h(H),h(H_0))\circ h.$$

Next, employing again DOI techniques, one proves that this definition of the SSF 
\emph{does not depend} on the function $h$ and, moreover,
\begin{itemize}\item[$(\alpha)$]
 if~$H\geq H_0$, then $\xi^{(4)}(\, \cdot \,; H,H_0) \geq 0$, 
\item[$(\beta)$] $\xi^{(4)}(\, \cdot \,; H,H_0) \in L_1(\bbR)$.
 \end{itemize}

\subsection{More general classes of perturbations.}
At this point we return to the case where the von Neumann algebra is the algebra $\cB(\cH)$ equipped with the standard trace and consider the situation when the perturbation is no longer a trace class operator. We note, that for the following results we \emph{will not specify} the class of functions $f$, for which the Krein trace formula holds. We are only interested in the  \emph{existence} of the SSF for a more general class of perturbations.

The first result, generalising the class of operators $H_0,H$ is due to M.G.Krein \cite{Kr62}.
\begin{theorem}[Resolvent comparable case]Let  the self-adjoint operators $H_0,H$ be such that
\begin{equation} 
[R_H(z)-R_{H_0}(z)] \in\cB_1(\cH), \quad z\in\rho(H_0)\cap\rho(H).   \lb{rc}
\end{equation} 
Then there exists a spectral shift function $\xi(\, \cdot \,; H,H_0)$, satisfying the weighted integrability condition
$$ \xi(\lambda;H,H_0) \in L^1\big(\bbR; (1+\lambda^2)^{-1} d\lambda\big).$$
\end{theorem}
We emphasize that in the present resolvent comparable case \eqref{rc}, this SSF is defined  only \emph{up to an additive constant}. 

Just as in the case of a trace class perturbation, the SSF for resolvent comparable operators $H_0,H$ possesses the following property:\\[1mm] 
$\bullet$ Suppose that in some interval $(a_0,b_0)$ the spectrum of $H_0$ is \emph{discrete} and let $\delta=(a,b)$, $a_0<a<b<b_0.$ Then the analogue of \eqref{xi_on_discrete} holds, that is,
\begin{equation}\label{xi_on_discrete_rc}
\xi(b_-; H,H_0) - \xi(a_+; H,H_0) = N_0(\delta)-N(\delta),
\end{equation}
where $N_0(\delta)$ (respectively, $N(\delta)$) are the sum of the multiplicities of the eigenvalues 
of $H_0$ (respectively, $H$) in $\delta$.

In the particular case of \emph{lower semibounded operators} $H_0$ and $H$ equality \eqref{xi_on_discrete_rc} allows us to \emph{naturally fix} the additive constant in the following way. To the left of the spectra of $H_0$ and $H$, the eigenvalue counting functions $N_0(\cdot)$ and $N(\cdot)$ are zero. Therefore, by equality \eqref{xi_on_discrete_rc} the SSF $\xi(\, \cdot \,; H,H_0)$ is a constant to the left of the spectra of $H_0$ and $H$, and it is custom to set this constant equal to zero,
$$
\xi(\lambda; H,H_0) = 0, \quad \lambda < \inf(\sigma(H_0) \cup \sigma(H)). 
$$

In the following we describe a particular way to introduce the SSF for the pair $(H,H_0)$ by what is usually called the \emph{invariance principle}. We note that this principle was used in a construction in \cite{PSZ} of the SSF for trace class perturbations at Step~$(iv)$, where we passed to unbounded operators.

Let $\Omega$ be an interval containing the spectra of $H_0$ and $H$, and let $\phi$ be an arbitrary bounded monotone  ``sufficiently" smooth function on $\Omega$. Suppose that
$$ [\phi(H)-\phi(H_0)] \in \cB_1(\cH)$$
then, the SSF $\xi(\, \cdot \,; H,H_0)$ can be defined as follows:
\begin{equation}\label{invar}
\xi(\lambda;H,H_0)=\sgn\big(\phi'(\lambda)\big)\xi(\phi(\lambda);\phi(H),\phi(H_0)).
\end{equation}
For the function $\xi(\, \cdot \,; H,H_0)$ the Lifshitz--Krein trace formula \eqref{trace_formula} holds for some class of admissible functions $f$. The latter class depends on $\phi$. 

We note the following result (see \cite[Sect.~8.11]{Ya92}):

\begin{proposition}Let $(H-H_0) \in\cB_1(\cH)$. Then the spectral shift functions for the pairs $(H,H_0)$ and $(\phi(H),\phi(H_0))$ are associated via equality \eqref{invar} up to an additive, integer-valued constant.
\end{proposition}

The methods of construction of the SSF introduced in this survey are only a sample of a plethora of possibilities. There are many others, which we did not cover here. We only mention a few of them:

\begin{itemize}
\item Sobolev \cite{So93} suggested a way of constructing the SSF for trace class perturbations via the ``argument of the perturbation determinant''. This construction allows one to establish \emph{pointwise} estimates on the SSF, and, in some cases, proves continuity of SSF on the absolutely continuous spectrum of $H_0$ (the latter coincides with that of $H$).
\item Koplienko \cite{Ko86} proved the existence of the SSF for the pair of operators 
$(H,H_0)$  under the assumption that for some $\varepsilon>0, 1\leq p<\infty$
$$[R_H(z)-R_{H_0}(z)](H_0^2+i)^{-\varepsilon}\in\cB_1(\cH), \quad 
[R_H(z)-R_{H_0}(z)] \in \cB_p(\cH).$$
\item Yafaev \cite{Ya05} proved that the SSF exists for a pair of operators $H_0$ and $H$ satisfying the assumption that for some $m \in \bbN$, $m$ odd, 
$$ \big[R_H^m(z)-R_{H_0}^m(z)\big] \in\cB_1(\cH).$$
\item Koplienko \cite{Ko84} proposed another function, which is called the Koplienko SSF \cite{Ko84}, and is constructed under the assumption that $(H-H_0)\in\cB_2(\cH)$ (see also 
\cite{GPS08}). For recent developments of this line of thought, we refer to \cite{PSS13}.
\end{itemize}

\section{The Witten index}

The Witten index of an operator $T$ in a complex separable Hilbert space $\cH$  provides a generalisation of  the 
Fredholm index of $T$ in certain cases where the operator $T$ ceases to have the 
Fredholm property. 
The Witten index possesses stability 
properties with respect to additive perturbations, which are analogous to, but more restrictive than,
the stability properties of the Fredholm index (roughly speaking, only relative trace class perturbations, as opposed to relative compact ones, are permitted). After the publication
of \cite{Wi82} this notion became popular in connection with a variety of 
examples in supersymmetric quantum theory in the 1980's. One reason the Witten index
has attracted little attention in recent years is that its connection with geometric questions remains
unclear (see however \cite{BMS88}, \cite{CK15}). This is a matter deserving further investigation. 
For more historical details we refer to the paragraphs following Theorem \ref{tW}. 

First, we recall the definitions and some of the basic properties of the Witten index. In the next section, we will derive new properties of 
the Witten index of a certain model operator.

We start with the following facts on trace class properties of resolvent and semigroup differences.

Then the following well-known and standard assertions hold for resolvent and semigroup comparable operators (see item $(ii)$ below): 

\begin{proposition}
Suppose that $0 \leq S_j$, $j=1,2$, are nonnegative, self-adjoint operators in $\cH$. \\[1mm] 
$(i)$ If $\big[(S_2 - z_0 )^{-1} - (S_1 - z_0 )^{-1}\big] \in \cB_1(\cH)$ for 
some $z_0 \in \rho(S_1)\cap \rho(S_2)$, then actually, 
\begin{equation*}
\big[(S_2 - z )^{-1} - (S_1 - z )^{-1}\big] \in \cB_1(\cH) 
\, \text{\emph{ for all} $z \in \rho(S_1)\cap \rho(S_2)$.}
\end{equation*} 
$(ii)$ If $\big[e^{- t_0 S_2} - e^{- t_0 S_1}\big] \in \cB_1(\cH)$ for some $t_0 > 0$, then 
\begin{equation*}
\big[e^{- t S_2} - e^{- t S_1}\big] \in \cB_1(\cH) \, \text{ \emph{for all} $t \geq t_0$.}  
\end{equation*}
\end{proposition}

The preceding fact allows one to consider the following two definitions.

Let $T$ be a closed, linear, densely defined operator in $\cH$. 
 Suppose that for some $($and hence for all\,$)$ 
$z \in \bbC \backslash [0,\infty) \subseteq [\rho(T^*T) \cap \rho(TT^*)]$,  
\begin{equation*} 
\big[(T^* T - z)^{-1} - (TT^* - z)^{-1}\big] \in \cB_1(\cH). 
\end{equation*}  
Then one introduces the resolvent regularization 
\begin{equation*}
\Delta_r(T, \lambda) = (- \lambda) \tr_{\cH}\big((T^* T - \lambda)^{-1}
- (T T^* - \lambda)^{-1}\big), \quad \lambda < 0.     
\end{equation*} 
\emph{The resolvent regularized Witten index $W_r (T)$} of $T$ is then defined by  
\begin{equation*} 
W_r(T) = \lim_{\lambda \uparrow 0} \Delta_r(T, \lambda), 
\end{equation*}
whenever this limit exists. 

Similarly, suppose that for some $t_0 > 0$  
\begin{equation*}
\big[e^{- t_0 T^* T} - e^{- t_0 TT^*}\big] \in \cB_1(\cH).   
\end{equation*} 
Then $\big(e^{-t T^*T} - e^{-t TT^*}\big) \in \cB_1(H)$ for all $t >t_0$ 
and one introduces the semigroup regularization 
\begin{equation*}
\Delta_s(T, t) = \tr_{\cH}\big(e^{-t T^*T} - e^{-t TT^*}\big), \quad t > 0.    
\end{equation*} 
\emph{The semigroup regularized Witten index $W_s (T)$} of $T$ is then defined by  
\begin{equation*} 
W_s(T) = \lim_{t \uparrow \infty} \Delta_s(T, t),   
\end{equation*}
whenever this limit exists.

One recalls that a densely defined and closed operator $T$ in a Hilbert space $\cH$ is said to be \emph{Fredholm} if 
$\ran(T)$ is closed and  $\dim(\ker(T)) + \dim(\ker(T^*)) < \infty$. In this case, the \emph{Fredholm index} $\ind(T):=\dim(\ker(T))-\dim(\ker(T^*))$.
The following result, obtained in \cite{BGGSS87} and  \cite{GS88}, states that both (resolvent and semigroup) regularized Witten indices coincide with the Fredholm index in the special  case of Fredholm operators. 

\begin{theorem}
Let $T$ be an $($unbounded\,$)$ Fredholm operator in $H$. Suppose that  
$\big[(T^* T - z )^{-1} - (TT^* - z )^{-1}\big],\
\big[e^{- t_0 T^* T} - e^{- t_0 TT^*}\big] \in \cB_1(\cH)$ for some $z \in \bbC \backslash [0,\infty),$ and $\ t_0>0$. Then 
$$ \ind(T)=W_r(T)=W_s(T).$$
\end{theorem}

In general (i.e., if $T$ is not Fredholm), $W_r(T)$ (respectively, $W_s(T)$) is not necessarily 
integer-valued; in fact, it can be any real number. As a concrete 
example, we mention the two-dimensional  magnetic field system discussed by 
Aharonov and Casher \cite{AC79} which demonstrates that the 
resolvent and semigroup regularized Witten indices have the meaning of 
(non-quantized) magnetic flux $F \in \bbR$ which indeed can be  any 
prescribed real number.   

Expressing the Witten index $W_s(T)$ (respectively, $W_r(T)$) of an operator $T$ in terms of the spectral shift function $\xi(\, \cdot \,; T^*T, TT^*)$ requires of course the choice of a concrete representative of the SSF:

\begin{theorem} \cite{BGGSS87,GS88} \lb{tW}
$(i)$ Suppose that $\big[e^{- t_0 T^* T} - e^{- t_0 TT^*}\big] \in \cB_1(\cH),\ t_0>0$ and the SSF $\xi(\, \cdot \,; T^*T,TT^*)$
, uniquely defined by the requirement $\xi(\lambda;T^*T,TT^*)=0,\ \lambda<0$,
 is continuous from above at $\lambda=0$. Then the semigroup regularized Witten index $W_s(T)$ of $T$ exists and 
$$W_s(T)=-\xi(0_+; T^*T,TT^*).$$  
$(ii)$ Suppose that $\big[(T^* T - z )^{-1} - (TT^* - z )^{-1}\big],\ z \in \bbC \backslash [0,\infty)$ and $\xi(\, \cdot \,; T^*T,TT^*)$, uniquely defined by the requirement $\xi(\lambda;T^*T,TT^*)=0,\ \lambda<0$,
 is bounded and piecewise continuous on $\bbR$. Then the resolvent regularized Witten index $W_s(T)$ of $T$ exists and 
$$W_r(T)=-\xi(0_+; T^*T,TT^*).$$ 
\end{theorem}

The first relations between index theory for not necessarily Fredholm operators and the 
Lifshitz--Krein  spectral shift function were established in \cite{BGGSS87}, \cite{Ge86}, \cite{GS88}, 
and independently in \cite{CP86}. In fact, inspired by index calculations of Callias \cite{Ca78} in connection with noncompact manifolds, the more general notion of the Witten index was 
studied and identified 
with the value of an appropriate spectral shift function at zero in \cite{BGGSS87} and 
\cite{GS88} (see also \cite{Ge86}, \cite[Ch.\ 5]{Th92}). Similiar investigations in search of an index theory for non-Fredholm operators were untertaken in \cite{CP86} in a slightly different direction,  based on principal functions and their connection to Krein's spectral shift function.  

The index calculations by Callias created considerable  interest, especially, in connection with certain aspects of supersymmetric quantum mechanics. Since a detailed list of references in this context is beyond the scope of this paper we only refer to \cite{An89}, \cite{An90}, \cite{An90a}, \cite{An93}, \cite{An94}, \cite{BS78}, \cite{BB84}, \cite{Bu95}, \cite{EGH80}, \cite{FOW87}, 
\cite{Hi83}, \cite{Hi86}, \cite{IM84}, \cite{Ko11}, \cite{Ko15}, \cite{NS86}, \cite{NS86a}, 
\cite[Ch.\ 5]{Th92} and the detailed lists of references cited therein. While \cite{BGGSS87} and \cite{Ge86} focused on index theorems for concrete one and 
two-dimensional  supersymmetric systems (in particular, the trace formula \eqref{ptf_eq} and the function $g_z(\cdot)$ in \eqref{gz} were discussed in 
\cite{BGGSS87} and \cite{Ge86} in the special case where $\cH = \bbC$), \cite{GS88} treated abstract Fredholm and Witten indices in terms of the spectral shift function and proved their invariance with respect to appropriate classes of perturbations. Soon after, a general abstract approach to supersymmetric scattering theory involving the spectral shift function was developed in \cite{BMS88} (see also \cite{Bu92}, \cite[Chs.\ IX, X]{Mu87}, \cite{Mu88}) and applied to relative index theorems in the context of manifolds Euclidean at infinity. 

\begin{example}\label{exampleBGGSS}
As an example of practical use of the abstract results, \cite{BGGSS87} considered the operator 
$$T=\frac{d}{dt}+M_\theta,\quad \dom(T)=W^{2,1}(\bbR),$$
acting on the standard Hilbert space $L^2(\bbR)$. Here $W^{2,1}(\bbR)$ is the Sobolev space, $M_\theta$ is the operator of multiplication by a bounded function $\theta$ on $\bbR$. Assuming existence of the asymptotes $\lim_{t\rightarrow\pm \infty}\theta(t)=\theta_{\pm}\in\bbR,$ 
and some additional conditions on $\theta$, it is shown in \cite{BGGSS87} that for the resolvent regularization one obtains 
\begin{align}\label{example_BGGSS} 
\begin{split} 
\Delta_r(T, \lambda) &= (- \lambda) \tr_{\cH}\big((T^* T - \lambda)^{-1}
- (T T^* - \lambda)^{-1}\big)   \\
&= \big[\theta_+(\theta_+^2-\lambda)^{-1/2}+\theta_-(\theta_-^2-\lambda)^{-1/2}\big]\big/2, 
\quad \lambda\in\bbC \backslash [0,\infty),
\end{split} 
\end{align}
and therefore, 
$$W_r(T)= [\sgn(\theta_+)-\sgn(\theta_-)]/2.$$
\end{example}

Next, we view  the operator $T$ from the preceding example as an operator of the form 
$T=\bsD_\bsA^{}=\frac{d}{dt}+\bsA$ on the Hilbert space $L^2(\bbR;\bbC)$, where $\bsA$ is the operator generated by the family of operators $\{A(t)\}_{t\in\bbR}$ on the Hilbert space $\bbC$, with $A(t)$ given by multiplication by $\theta(t)$, $t\in\bbR$.

Our main objective in this survey is to consider a more general situation where the family $\{A(t)\}$ consists of  operators acting on an arbitrary complex, separable initial Hilbert space $\cH$, and the resulting operator $\bsD_\bsA^{}=\frac{d}{dt}+\bsA$ acts on the Hilbert space $L^2(\bbR;\cH)$.
Operators of this form $\bsD_\bsA^{}$ arise in connection
with Dirac-type operators (on compact and noncompact manifolds), the Maslov
index, Morse theory (index), Floer homology, winding numbers, Sturm's oscillation
theory, dynamical systems, etc. (cf.\ \cite{GLMST11} and the extensive list of references therein).

To date, strong conditions need to be imposed on the family $A(t)$ in order to obtain
 the resolvent and semigroup Witten indices of $\bsD_\bsA^{}$ and express them in terms of the spectral shift function for the asymptotes $A_\pm$ of the family $A(t)$ as $t\to \pm \infty$. 
The following is the main hypothesis, under which the results stated below are proved.

\begin{hypothesis}\label{hyp} 
$(i)$ Suppose $\cH$ is a complex, separable Hilbert space. \\[1mm] 
$(ii)$ Assume $A_-$ is a self-adjoint operator on $\dom(A_-) \subseteq \cH$. \\[1mm] 
$(iii)$ Suppose there exists a family of bounded self-adjoint operators $B(t)$, $t\in\bbR$ with $t\mapsto B(t)$ weakly locally absolutely continuous on $\bbR$, implying the 
existence of a family of bounded self-adjoint operators $\{B'(t)\}_{t\in\bbR}$ in $\cH$ such that for 
a.e.\ $t\in\bbR$,  
\begin{equation*} 
\frac{d}{dt} (g,B(t) h)_{\cH} = (g,B'(t) h)_{\cH}, \quad g, h\in\cH.  
\end{equation*} 
$(iv)$ Assume that the family $\{B'(t)\}$ is \emph{$A_-$-relative trace class}, that is, suppose that $B'(t)(|A_-| + 1)^{-1}\in\cB_1(\cH)$, $t\in\bbR$. In addition, we assume that
\begin{equation*}  
\int_\bbR dt \, \big\|B'(t) (|A_-| + 1)^{-1}\big\|_{\cB_1(\cH)} < \infty.
\end{equation*}
\end{hypothesis}

\begin{remark}
$(i)$ We note that, in fact, the subsequent results hold in a more general situation, when the operators $B(t),\ t\in\bbR,$ are allowed to be unbounded and some additional measurability conditions of the families $\{B(t)\},$ \, $\{B'(t)\}$ are imposed. \\[1mm] 
$(ii)$ The assumption $(iv)$ above, that the operators $B'(t), \ t\in\bbR,$ are relative trace class, namely, $B'(t)(|A_-| + 1)^{-1}\in\cB_1(\cH)$, is the main assumption, which implies various trace relations below.  In Section \ref{dim_1} we will discuss an example where we replace the relative
trace class hypotheses with a  relative Hilbert--Schmidt class assumption. \hfill $\diamond$ 
\end{remark}

From this point on we always assume Hypothesis \ref{hyp}.

\subsection{Definition of the operator $\bsD_\bsA^{}$.} 
We introduce the family of self-adjoint operators 
$\{A(t)\}_{t\in\bbR}$ in $\cH$ by 
\begin{equation*}
A(t) = A_- + B(t), \quad \dom(A(t)) = \dom(A_-), \; t\in\bbR.
\end{equation*}
Hypothesis \ref{hyp} ensures the existence of the \emph{asymptote $A_+$} as $t\rightarrow\infty$ in the norm-resolvent sense, $\dom(A_+) = \dom(A_-)$, with $A_+$ self-adjoint in  $\cH$, that is 
$$\lim_{t\rightarrow+\infty}\big\|(A(t)-z)^{-1}-(A_+-z)^{-1}\big\|_{\cB(H)}=0.$$

Let $\bsA$ in $L^2(\bbR;\cH)$ be the operator associated with the family 
$\{A(t)\}_{t\in\bbR}$ in $\cH$ by
\begin{align*}
&(\bsA f)(t) = A(t) f(t) \, \text{ for a.e.\ $t\in\bbR$,}   \no \\
& f \in \dom(\bsA) = \bigg\{g \in L^2(\bbR;\cH) \,\bigg|\,
g(t)\in \dom(A(t)) \text{ for a.e.\ } t\in\bbR,     \\
& \quad t \mapsto A(t)g(t) \text{ is (weakly) measurable,} \,  
\int_{\bbR} dt \, \|A(t) g(t)\|_{\cH}^2 < \infty\bigg\}.    \no
\end{align*}

We define also the operator ${\bf d/dt}$ in $L^2(\bbR;\cH)$  by setting
\begin{align*}
\bigg({\bf\f{d}{dt}}f\bigg)(t) = f'(t) \, \text{ a.e.\ $t\in\bbR$,} 
\quad f\in \dom({\bf d/dt}) = W^{2,1}(\bbR;\cH),
\end{align*}
where $W^{1,2}(\bbR;\cH)$ denotes the usual Sobolev space of $L^2(\bbR;\cH)$-functions with the first distributional derivative in $L^2(\bbR;\cH)$. 

Now, we introduce the operator $\bsD_\bsA^{}$ in 
$L^2(\bbR;\cH)$ by setting
\begin{equation}\label{def_D_A}
\bsD_\bsA^{} = {\bf\f{d}{dt}} + \bsA,
\quad \dom(\bsD_\bsA^{})= \dom({\bf d/dt}) \cap \dom(\bsA).  
\end{equation}

\begin{proposition}\cite{GLMST11}
Assume Hypothesis \ref{hyp}. Then the operator
$\bsD_\bsA^{}$ is densely defined and closed in $L^2(\bbR; \cH)$ and 
the adjoint $\bsD_\bsA^*$ of $\bsD_\bsA^{}$ in $L^2(\bbR; \cH)$ is given 
by
$$\bsD_\bsA^*=- {\bf\f{d}{dt}} + \bsA, \quad
\dom(\bsD_\bsA^*) = \dom(\bsD_\bsA^{}).$$
\end{proposition}

\subsection{The  principle trace formula.}
The following result relates the trace of the difference of the 
resolvents of positive operators $|\bsD_\bsA^{}|^2$ and $|\bsD_\bsA^{*}|^2$ in $L^2(\bbR;\cH)$, and the trace of the difference of $g_z(A_+)$ and $g_z(A_-)$ in $\cH$, where
\begin{equation} 
g_z(x) = x(x^2-z)^{-1/2}, \quad x\in\bbR, \; z\in\C\backslash [0,\infty).    \lb{gz} 
\end{equation}

\begin{theorem}\cite{GLMST11}\label{ptf} Assume Hypothesis \ref{hyp}. Then, 
\begin{align*}
& \big[\big(|\bsD_\bsA^{*}|^2 - z\big)^{-1}-\big(|\bsD_\bsA^{}|^2 - z\big)^{-1}\big] \in \cB_1\big(L^2(\bbR;\cH)\big),  
\quad z\in\rho(|\bsD_\bsA^{}|^2) \cap \rho(|\bsD_\bsA^{*}|^2),  \\
& [g_z(A_+)-g_z(A_-)] \in \cB_1(\cH),  \quad 
z\in\rho\big(A_+^2\big) \cap \rho\big(A_-^2\big),
\end{align*}
and the following principal trace formula holds for all $z\in\bbC\backslash [0,\infty)$,  
\begin{align}\label{ptf_eq}
   \tr_{L^2(\bbR;\cH)}\big(\big(|\bsD_\bsA^{*}|^2 - z\big)^{-1}-\big(|\bsD_\bsA^{}|^2 - z\big)^{-1}\big) = \frac{1}{2z} \tr_\cH \big(g_z(A_+)-g_z(A_-)\big).
\end{align} 
\end{theorem}

\begin{remark}
$(i)$ Pushnitski \cite{Pu08} was the first to investigate, under the more restrictive assumption 
that the operators $B(t)$ are \emph{trace class}, 
a trace formula of this kind. In our more general setting of relative trace class perturbations, this formula is obtained in \cite{GLMST11} by 
an approximation technique on both sides of the equation and a non-trivial DOI technique. \\[1mm] 
$(ii)$ Employing basic notions in scattering theory and the Jost--Pais-type reduction of Fredholm determinant, a recent paper \cite{CGPST14} provides a new proof of the principle trace formula in the case of trace class perturbations. \\[1mm] 
$(iii)$ If $\cH=\bbC$, the principal trace formula yields \eqref{example_BGGSS} for the example considered by D.Bolle {\it et al} for $\phi_\pm=\pm1$. \hfill $\diamond$ 
\end{remark}

\subsection{Pushnitki's formula relating two SSFs.}
Having at hand the principal trace formula, we now aim at correlating the underlying SSFs, 
$\xi(\, \cdot \,; |\bsD_\bsA^{*}|^2,|\bsD_\bsA^{}|^2) $  and  $\xi(\, \cdot \,; A_+,A_-).$

First, we need to properly introduce the SSF 
$\xi\big(\, \cdot \,; |\bsD_\bsA^{*}|^2, |\bsD_\bsA^{}|^2\big)$ associated with the pair of positive operators $|\bsD_\bsA^{*}|^2$ and $|\bsD_\bsA^{}|^2$.
By Theorem \ref{ptf} we have 
 \begin{equation*}
\Big[\big(|\bsD_\bsA^{*}|^2 - z\big)^{-1}-\big(|\bsD_\bsA^{}|^2 - z\big)^{-1}\Big] \in \cB_1 
\big(L^2(\bbR;\cH)\big),
\end{equation*}
and therefore, one can uniquely introduce 
$\xi\big(\, \cdot \,; |\bsD_\bsA^{*}|^2, |\bsD_\bsA^{}|^2\big)$ by requiring that
\begin{equation*}
\xi\big(\lambda; |\bsD_\bsA^{*}|^2, |\bsD_\bsA^{}|^2\big)= 0, \quad \lambda < 0,    \lb{2.27}
\end{equation*}
and then obtains
\begin{equation*}
\tr\Big(\big(|\bsD_\bsA^{*}|^2 - z\big)^{-1}-\big(|\bsD_\bsA^{}|^2 - z\big)^{-1}\Big)
= - \int_{[0, \infty)}  \frac{\xi\big(\lambda; |\bsD_\bsA^{*}|^2, |\bsD_\bsA^{}|^2\big) \, 
d\lambda}{(\lambda -z)^2},
\end{equation*}
for all $z\in\bbC\backslash [0,\infty).$

We shall introduce the spectral shift function associated with the pair 
$(A_+, A_-)$ via the invariance principle. 
By Theorem \ref{ptf}, $[g_{-1}(A_+)-g_{-1}(A_-)] \in \cB_1(\cH)$ 
and so one can define the SSF $\xi(\, \cdot \,; A_+,A_-)$ for the pair $(A_+, A_-)$ by setting
\begin{equation*}
\xi(\nu; A_+, A_-) := \xi(g_{-1} (\nu); g_{-1}(A_+), 
g_{-1}(A_-)), \quad \nu\in\bbR, 
\end{equation*}
 where
$\xi(\,\cdot\,;g_{-1}(A_+), g_{-1}(A_-))$ is the spectral shift function associated with the pair
$(g_{-1}(A_+), g_{-1}(A_-))$ uniquely determined by the requirement 
\begin{equation*}
\xi(\,\cdot\,;g_{-1}(A_+), g_{-1}(A_-)) \in L^1(\bbR; d\omega). 
\end{equation*} 
Therefore, by the Lifshitz--Krein trace formula \eqref{trace_formula},  
\begin{equation*}
\tr_{\cH}\big(g_{z}(A_+) - g_{z}(A_-)\big)
  = - z \int_{\bbR} \frac{\xi(\nu; A_+, A_-) \, d\nu}{(\nu^2 - z)^{3/2}},
\quad  z\in\bbC\backslash [0,\infty).     
\end{equation*}

By the principal trace formula one obtains the equality 
\begin{align*}
\int_{[0, \infty)}  &\frac{\xi\big(\lambda; |\bsD_\bsA^{*}|^2, |\bsD_\bsA^{}|^2\big) \, 
d\lambda}{(\lambda -z)^{-2}}\\
&=-\tr_{L^2(\bbR;\cH)} \Big(\big(|\bsD_\bsA^{*}|^2 - z \, \bsI\big)^{-1} - \big(|\bsD_\bsA^{}|^2 - z \, 
\bsI\big)^{-1}\Big)\\
&=-\frac1{2z}\tr_{\cH}\big(g_{z}(A_+) - g_{z}(A_-)\big)\\
&= \frac{1}{2} \int_{\bbR} \frac{\xi(\nu; A_+, A_-) \, d\nu}{(\nu^2 - z)^{3/2}},
\quad z\in\bbC\backslash [0,\infty),
\end{align*}
or, equivalently, 
\begin{align*}
 & \int_{[0, \infty)} \xi\big(\lambda; |\bsD_\bsA^{*}|^2, |\bsD_\bsA^{}|^2\big) 
\bigg(\frac{d}{dz}(\lambda -z)^{-1}\bigg)
  d\lambda\\
& \quad = \int_{\bbR} \xi(\nu; A_+, A_-)\bigg(\frac{d}{dz} (\nu^2 - z)^{-1/2}\bigg) d\nu.   
\end{align*} 
Integrating the preceding equality with respect to $z$ from a fixed point $z_0 
\in (-\infty,0)$ to $z\in\bbC\backslash\bbR$, along a straight line connecting $z_0$ and 
$z$, yields
\begin{align*}
& \int_{[0, \infty)} \xi\big(\lambda; |\bsD_\bsA^{*}|^2, |\bsD_\bsA^{}|^2\big)
\bigg(\frac{1}{\lambda - z} - \frac{1}{\lambda - z_0}\bigg) d\lambda      \\
& \quad = \int_{\bbR} \xi(\nu; A_+, A_-)\big[(\nu^2 - z)^{-1/2} - (\nu^2 - 
z_0)^{-1/2}\big] \, d\nu,
\quad z\in\bbC\backslash [0,\infty).  
\end{align*}
Applying the Stieltjes inversion formula then permits one to express the SSF function 
$\xi\big(\, \cdot \,; |\bsD_\bsA^{*}|^2, |\bsD_\bsA^{}|^2\big)$ in terms of 
$\xi(\, \cdot \,; A_+, A_-)$ as follows, 
\begin{align*}
\xi\big(\lambda; |\bsD_\bsA^{*}|^2, |\bsD_\bsA^{}|^2\big)
&= \lim_{\varepsilon\downarrow 0} \frac{1}{\pi} \int_{[0,\infty)}
\xi\big(\lambda'; |\bsD_\bsA^{*}|^2, |\bsD_\bsA^{}|^2\big) \Im\big(((\lambda'  - 
\lambda) - i \varepsilon)^{-1}\big) d\lambda'
\no  \\
& =  \lim_{\varepsilon\downarrow 0} \frac{1}{\pi}  \int_{\bbR} \xi(\nu; A_+, A_-)
\Im\big((\nu^2 - \lambda - i \varepsilon)^{-1/2}\big) d\nu     \no  \\
& = \frac{1}{\pi}  \int_{- \lambda^{1/2}}^{\lambda^{1/2}}
\frac{\xi(\nu; A_+, A_-) \, d\nu}{(\lambda - \nu^2)^{1/2}} \, \text{ for 
a.e.\ $\lambda > 0$.}  
\end{align*}
In the last equality here one should be careful with various estimates in order to apply Lebesgue's dominated convergence theorem. We omit further details and refer to \cite{GLMST11}. 

Putting all of this together we have the following remarkable formula, which expresses the SSF, 
$\xi\big(\, \cdot \,; |\bsD_\bsA^{*}|^2, |\bsD_\bsA^{}|^2\big)$, in terms of the SSF 
$\xi(\, \cdot \,; ; A_+,A_-)$. It is this formula that allows us to express (Fredholm/Witten) index of the operator $\bsD_\bsA^{}$ in terms of the spectral shift function $\xi(\, \cdot \,;  A_+,A_-)$.
Note, that this formula can be viewed as an Abel-type transform. 

\begin{theorem}\cite{GLMST11}\label{thm_Push}
Assume Hypothesis \ref{hyp} and define the spectral shift functions  
$\xi\big(\, \cdot \,;  |\bsD_\bsA^{*}|^2, |\bsD_\bsA^{}|^2\big)$ and 
$\xi(\, \cdot \,; ; A_+,A_-)$ as above. Then,  for a.e.\ $\lambda>0$,
\begin{equation*} 
\xi\big(\lambda; |\bsD_\bsA^{*}|^2, |\bsD_\bsA^{}|^2\big)=\frac{1}{\pi}\int_{-\lambda^{1/2}}^{\lambda^{1/2}}
\frac{\xi(\nu; A_+,A_-) \, d\nu}{(\lambda-\nu^2)^{1/2}}, 
\end{equation*} 
with a convergent Lebesgue integral on the right-hand side. 
\end{theorem}

A  formula of this kind was originally obtained for trace class perturbations $B(t)$ by 
Pushnitski \cite{Pu08}.

\subsection{The Fredholm case.} 
In order to study the Witten index of the operator  $\bsD_\bsA^{}$ we first need to understand under which additional assumptions this operator is Fredholm, which is of course the simpler case. The following result yields necessary and sufficient conditions for the latter.

\begin{theorem}\cite[Theorem 2.6]{CGPST15}\label{thn_Fredholm} Assume Hypothesis 
\ref{hyp}. Then the  operator  $\bsD_\bsA^{}$ is Fredholm if and only if 
$0 \in \rho(A_+)\cap\rho(A_-)$ $($i.e., $A_\pm$ are both boundedly invertible\,$)$. 
\end{theorem}

In fact, this theorem yields a complete description of the essential
spectrum of $\bsD_\bsA^{}$.  Here we define the essential spectrum of a
densely defined, closed, linear operator $T$ in a complex, separable Hilbert space $\cH$ as
\begin{equation*}
\sigma_{ess}(T) = \{\lambda \in \bbC \, | \, (T-\lambda I_{\cH})  
\text{ is not Fredholm}\}
\end{equation*}
(but caution the reader that several inequivalent, yet meaningful, definitions of essential spectra for non-self-adjoint operators exist, see, e.g., \cite[Ch.~IX]{EE89}). 

\begin{corollary}\cite[Corollary 2.8]{CGPST15} Assume Hypothesis \ref{hyp}. Then, 
\label{spec_ess} 
$$\sigma_{{\rm ess}}(\bsD_\bsA^{})= (\sigma(A_+) + i \, \bbR) \cup 
(\sigma(A_-) + i \, \bbR).
$$
\end{corollary}

By Theorem \ref{thn_Fredholm}, when the operator $\bsD_\bsA^{}$  is Fredholm, we have that $0 \in \rho(A_+)\cap\rho(A_-)$. Thus, by Corollary \ref{spec_ess}, $|\bsD_\bsA^{}|^2$ and 
$|\bsD_\bsA^{*}|^2$ have a gap in their essential spectrum near zero, that is, there exists an 
$a>0$ such that 
$$\sigma_{\rm ess}\big(|\bsD_\bsA^{}|^2\big)=\sigma_{\rm ess}\big(|\bsD_\bsA^{*}|^2\big)\subset [a,\infty).$$
This means that, on the interval $[0,a)$, the operators $|\bsD_\bsA^{}|^2$ and $|\bsD_\bsA^{*}|^2$ have discrete spectra. 
Hence, using  properties of the spectral shift function for discrete spectra (see property $(iii)$ in Subsection \ref{sec_properties_ssf}) one infers that 
\begin{equation*}
\xi\big(\lambda; |\bsD_\bsA^{*}|^2, |\bsD_\bsA^{}|^2\big) 
= \xi(0_+; |\bsD_\bsA^{*}|^2, |\bsD_\bsA^{}|^2), \quad
\lambda \in (0,\lambda_0),  
\end{equation*}
for $\lambda_0 < \inf(\sigma_{\rm ess} (|\bsD_\bsA^{}|^2)) = \inf(\sigma_{\rm ess} (|\bsD_\bsA^{*}|^2))$. 

On the other hand, since  $0 \in \rho(A_+)\cap\rho(A_-)$, there exists a constant 
$c\in\bbR$ such that $\xi(\,\cdot\,; A_+, A_-) = c$ a.e.\ on the interval 
 $(-\nu_0,\nu_0)$ for $0 < \nu_0$ sufficiently small (see property $(i)$ in Subsection \ref{sec_properties_ssf}).  Hence, the value $\xi(0; A_+, A_-)$ is well defined and 
\begin{equation*}
\xi(\nu; A_+, A_-) = \xi(0; A_+, A_-), \quad \nu \in (-\nu_0,\nu_0),
\end{equation*}
in particular, $\lim_{\nu\rightarrow 0}\xi(\nu; A_+,A_-)=\xi(0; A_+,A_-)$.

Thus, taking $\lambda \downarrow 0$ in Pushnitski's formula one infers  
\begin{align*}
\xi(0_+;& |\bsD_\bsA^{*}|^2, |\bsD_\bsA^{}|^2)=\lim_{\lambda\downarrow 0}\xi(\lambda; |\bsD_\bsA^{*}|^2, |\bsD_\bsA^{}|^2)\\
&=
\lim_{\lambda\downarrow 0+}\frac{1}{\pi}\int_{-\lambda^{1/2}}^{\lambda^{1/2}}
\frac{\xi(\nu; A_+,A_-) \, d\nu}{(\lambda-\nu^2)^{1/2}}=\lim_{\nu\rightarrow 0}\xi(\nu; A_+,A_-)\\
&=\xi(0; A_+,A_-)
\end{align*}
since $\pi^{-1}\int_{-\lambda^{1/2}}^{\lambda^{1/2}}
d\nu \, (\lambda- \nu^2)^{-1/2} = 1 $ for all $\lambda > 0$.

Thus, we obtain the following result linking the Fredholm index for $\bsD_\bsA^{}$ and the value of the SSF $\xi(\, \cdot \,; A_+, A_-)$ at zero.

\begin{theorem}\cite{GLMST11}
\label{thm_Fredholm_case}Assume Hypothesis \ref{hyp} and introduce the SSFs 
$\xi(\,\cdot\,; A_+, A_-)$ and $\xi\big(\, \cdot \,; |\bsD_\bsA^{*}|^2, |\bsD_\bsA^{}|^2\big)$ as above. Moreover, suppose that $0 \in 
\rho(A_+)\cap\rho(A_-)$. Then $\bsD_\bsA^{}$
is a Fredholm operator in $L^2(\bbR;\cH)$ and
\begin{equation}\label{eq_for_FI}
W_r(\bsD_\bsA^{})=\ind(\bsD_\bsA^{})= \xi\big(0_+; |\bsD_\bsA^{*}|^2, |\bsD_\bsA^{}|^2\big) 
= \xi(0; A_+, A_-).
\end{equation}
\end{theorem}

We emphasize that the assumption $0 \in 
\rho(A_+)\cap\rho(A_-)$ is crucial in the Fredholm index formula \eqref{eq_for_FI} of the operator $\bsD_\bsA^{}$. This assumption allows us to \emph{define} the value of SSF $\xi(\, \cdot \,; A_+, A_-)$ at zero. Generally speaking, the SSF $\xi(\, \cdot \,; A_+, A_-)$ is defined as an element in $L^1\big(\bbR; (|\nu|+1)^{-3}\big)$ (the space of classes of functions), so it does not make sense to speak of its value at a fixed point.

\subsection{Connection to spectral flow.}

The  relationship between spectral flow and the Fredholm index was first raised in the original articles of Atiyah--Patodi--Singer \cite{APS76}.
A definitive treatment of the question  for certain families of self-adjoint unbounded operators with compact resolvent was provided in \cite{RS95}, essentially, using the model operator formalism that we described above. For 
partial differential operators on noncompact manifolds it is typically the case that they possess some essential spectrum so that \cite{RS95} is not applicable. 
This motivated the investigation in \cite{Pu08} and \cite{GLMST11}.  The first of these papers introduces new
methods and ideas, relating the index/spectral flow connection to scattering theory and the spectral shift function.
However the conditions imposed in \cite{Pu08} are too restrictive to allow a wide class of examples. 
New tools were introduced in \cite{GLMST11}  as is explained above.  A more detailed history of these matters may also be found in 
 \cite{CGPST14} which also contains results on an index theory for certain non-Fredholm operators using the model operator formalism above.

One of the 
principle aims of \cite{GLMST11} was to extend the  results in \cite{RS95}
(albeit subject to  a relative trace class perturbation condition), in a fashion permitting essential spectra. This has motivated our interest
in the problem of applying these new methods to Dirac-type operators
on non-compact manifolds. There is a difficulty, however, in that the relative trace class perturbation assumption 
is not satisfied in this context (even in the one-dimensional case). In the last section of this review we will address
this difficulty via a class of examples.

Spectral flow is usually discussed in terms of measuring the nett number of
eigenvalues of a one parameter family of Fredholm operators that
change sign as one moves along the path.
In fact we need a more precise definition and use the one due to
Phillips \cite{Ph96}

Consider a norm continuous path $F_t$, $t\in [0; 1],$ of bounded self-adjoint
Fredholm  operators  joining $F_1$ and $F_0$.
For each $t$, we let $P_t$ be the spectral projection of $F_t$ corresponding to
the non-negative reals. Then we can write $F_t = (2P_t -1)|F_t|$.
Phillips showed that if one subdivides the path into small intervals
$[t_j, t_{j+1}]$ such that $P_{t_j}$ and $P_{t_{j+1}}$ are ``close''  in the Calkin algebra, then they
form a Fredholm pair (i.e., $P_{t_{j}}P_{t_{j+1}}$ is a Fredholm operator from $\ran (P_{t_{j+1}})$ to $\ran (P_{t_{j}})$)  and the spectral flow along $\{F_t\}_{t \in [0,1]}$ is defined by 
$$\sum_{j}\ind(P_{t_j}P_{t_{j+1}} : \ran (P_{t_{j+1}}) \to \ran (P_{t_j}) )$$

We will now state the main result of \cite{GLMST11} on the connection between  the spectral flow along the path $\{A(t)\}_{t \in \bbR}$ with the spectral shift functions and the Fredholm index of the model operator $\bsD_\bsA^{}$ introduced previously. However, first we need some preparatory observations.

 First, we note that spectral flow along the path of unbounded operators $\{A(t)\}_{t \in \bbR}$
is defined in terms of the flow along the bounded transforms 
$\big\{F_t=A(t) (I+A(t)^2)^{-1/2}\big\}_{t \in \bbR}$ using Phillips' definition above. Next, we remark that the fact that the spectral shift function is relevant to the discussion of spectral flow 
was first noticed by M\"uller \cite{Mu98} and explained in a systematic fashion  
in 2007 in \cite{ACS07}. There it was shown that, under very general conditions 
guaranteeing that both are well-defined, the spectral shift function and spectral flow
 are the same notion. The main technical tool exploited there 
is the theory of double operator integrals.

The new ingredient in \cite{GLMST11}  is a  formula, which connects the spectral flow with both the spectral shift function and the Fredholm index. This formula applies independently of whether the operators in the path have non-trivial essential spectrum or not.
More precisely, the spectral flow  along the family of Fredholm operators $\{A(t)\}_{t \in \bbR}$ coincides with the value of the Fredholm index of the operator $\bsD_\bsA^{}$ and the value of the SSF, $\xi(\, \cdot \,; A_+,A_-)$, computed at zero. 

\begin{theorem}\cite{GLMST11}
Assume Hypothesis \ref{hyp} and suppose that $0 \in \rho(A_+)\cap\rho(A_-)$. 
Then  $\big(E_{A_+}((-\infty,0)),E_{A_-}((-\infty,0))\big)$ form a Fredholm pair
and  the following equalities hold:
\begin{align*}
\text{\rm SpFlow} (\{A(t)\}_{t=-\infty}^\infty) & =\ind(E_{A_-}((-\infty,0)),E_{A_+}((-\infty,0)))\\
&={\tr}_{\cH}(E_{A_-}((-\infty,0))-E_{A_+}((-\infty,0))) \\
&=\xi(0; A_+,A_-) = \xi(0_+; \bsH_2, \bsH_1)    \lb{spf3A} =\ind (\bsD_\bsA^{}).  
\end{align*}
\end{theorem}

\section{Witten index: new results.} 

In the preceding section we discussed the notion of the Witten index and its connection with the Fredholm index as well as the spectral shift function. 
As we already know from Theorem \ref{thn_Fredholm}, the operator $\bsD_\bsA^{}$ is Fredholm if and only if $0\in \rho(A_+)\cap \rho(A_-)$, that is, the operators $A_\pm$ are both boundedly invertible. Moreover, if $0\in \rho(A_+)\cap \rho(A_-)$, then 
the Fredholm index can be computed as  
$$\ind(\bsD_\bsA^{})=\xi(0; A_+, A_-).$$

Here, the assumption $0 \in \rho(A_+)\cap\rho(A_-)$ is crucial, since in this case, there exists 
$0< \nu\in\bbR$, such that $\xi(\, \cdot \,; A_+,A_-)$ is constant  on the interval $(-\nu,\nu)$ so that one can speak about the value of the SSF, $\xi(\, \cdot \,; A_+, A_-)$, at zero. 
An important question then is to study an extension of index theory for the operator $\bsD_\bsA^{}$, when the latter ceases to be Fredholm. In this case $0\in\sigma(A_+),$ or $0\in\sigma(A_-)$ and therefore, the SSF $\xi(\, \cdot \,;A_+,A_-)$ \emph{is not  constant}, in general, on any interval 
$(-\nu,\nu), \nu>0$.

An approach to computing the Witten indices $W_{r}(\bsD_\bsA^{})$ (respectively, $W_{s}(\bsD_\bsA^{})$)  suggested in \cite{CGPST15} relies on the usage of right and left Lebesgue points of spectral shift functions. We start by briefly recalling this notion.

\begin{definition} \label{def_1}
Let $f \in L^1_{\loc} (\bbR; dx)$ and $h > 0$. \\
$(i)$ The point  $x \in \bbR$ is called a {\bf right Lebesgue point} of $f$ if 
there exists an $\alpha_+ \in \bbC$ such that  
\begin{equation*} 
\lim_{h \downarrow 0} \f{1}{h} \int_{x}^{x + h} dy \, |f(y) - \alpha_+| = 0.   
\end{equation*} 
We write $\alpha_+ = \Lf(x_+)$. \\
$(ii)$ The point $x \in \bbR$ is called a {\bf left Lebesgue point} of $f$ if 
there exists an $\alpha_- \in \bbC$ such that  
\begin{equation*} 
\lim_{h \downarrow 0} \f{1}{h} \int_{x - h}^{x} dy \, |f(y) - \alpha_-| = 0.   
\end{equation*}
We write $\alpha_- = \Lf(x_-)$. \\ 
$(iii)$ The point $x \in \bbR$ is called a {\bf Lebesgue point} of $f$ if 
there exist $\alpha \in \bbC$ such that  
\begin{equation*} 
\lim_{h \downarrow 0} \f{1}{2h} \int_{x - h}^{x + h} dy \, |f(y) - \alpha| = 0. 
\end{equation*} 
We write $\alpha = \Lf(x)$. That is, $x \in \bbR$ is a Lebesgue point 
of $f$ if and only if it is a left and a right Lebesgue point and 
$\alpha_+ = \alpha_- = \alpha$. 
\end{definition} 

We note that this definition of a Lebesgue point of $f$ is not universally adopted.  
For instance, \cite[p.~278]{HS65} define $x_0$ to be a Lebesgue point 
of $f$ if 
\begin{equation}\label{5.3}
\lim_{h \downarrow 0} \f{1}{h} \int_{0}^h dy \, |f(x + y) + f(x_0 - y) - 2 f(x)| = 0.   
\end{equation}
The elementary example 
\begin{equation*}
f(x;\beta) = \begin{cases} 0, & x<0, \\ \beta, & x=0, \\ 1, & x>0, \end{cases} 
\quad \beta \in \bbC,    
\end{equation*}
shows that 
$\Lf (0_+ ; \beta) = 1$ and $\Lf(0_- ; \beta) = 0$, that is, $x_0 = 0$ is a right and a left Lebesgue 
point of $f(\, \cdot \, ; \beta)$ in the sense of Definitions \ref{def_1}, whereas 
there exists no $\beta \in \bbC$ such that $f(\, \cdot \, ; \beta)$ satisfies \eqref{5.3} for $x_0 = 0$. 

Everywhere below we use the terms left and right Lebesgue point of a function in the sense of Definition \ref{def_1}.

\subsection{Connection between Lebesgue points of the SSFs, $\xi(\, \cdot \,; A_+,A_-)$ and 
$\xi\big(\, \cdot \,; |\bsD_\bsA^{*}|^2, |\bsD_\bsA^{}|^2\big)$.} 
As in the Fredholm case, the main ingredient in computing  the Witten index is Pushnitski's formula (see Theorem \ref{thm_Push}):

$$\xi\big(\lambda; |\bsD_\bsA^{*}|^2, |\bsD_\bsA^{}|^2\big)=\frac{1}{\pi}\int_{-\lambda^{1/2}}^{\lambda^{1/2}}
\frac{\xi(\nu; A_+,A_-) \, d\nu}{(\lambda-\nu^2)^{1/2}}.$$

We can rewrite this formula as follows:
$$\xi\big(\lambda; |\bsD_\bsA^{*}|^2, |\bsD_\bsA^{}|^2\big) = \frac{1}{\pi}\int_0^{\lambda^{1/2}}
\frac{d \nu \, [\xi(\nu; A_+,A_-) + \xi(-\nu; A_+,A_-)]}{(\lambda-\nu^2)^{1/2}},  
\quad \lambda > 0,    $$

and consider the operator 
 $S$, defined  by setting 
\begin{align} S: \begin{cases}
L^1_{loc}(\bbR; d\nu) \to L^1_{loc}((0,\infty); d\lambda),    \\[1mm] 
f \mapsto \frac1{\pi}\int_0^{\lambda^{1/2}} d \nu \, (\lambda-\nu^2)^{-1/2} f(\nu), \quad \lambda>0,
\end{cases}  
\end{align}
We then have the following result for the operator $S$:

\begin{lemma}\cite[Lemma 4.1]{CGPST15}\label{lem_3}
If $0$ is a right Lebesgue point for $f \in L^1_{loc}(\bbR; d\nu)$, then it is also a right Lebesgue 
point for $S f$ and $
\LSf (0_+) = \frac12\Lf (0_+).$
\end{lemma}

Hence, assuming that $0$ is a right and a left Lebesgue point of $\xi(,\cdot; A_+, A_-)$, an application of this lemma to the particular function $$f(\nu) = \xi(\nu,A_+,A_-)+\xi(-\nu,A_+,A_-),$$ 
$\nu > 0$, yields that $0$ is a right Lebesgue point of 
$\xi\big(\, \cdot \, ; |\bsD_\bsA^{*}|^2, |\bsD_\bsA^{}|^2\big)$ and  
\begin{equation}\label{eq_in_Lebesgue_sense}
\Lxi\big(0_+; |\bsD_\bsA^{*}|^2, |\bsD_\bsA^{}|^2\big)\stackrel{Lemma \, \ref{lem_3}}{=} 
[\Lxi(0_+; A_+,A_-) + \Lxi(0_-; A_+, A_-)]/2.
\end{equation}
Thus, in the case, when $0 \in\sigma(A_+)$ (or $\sigma(A_-)$), we can still correlate the values at zero of the functions $\xi(,\cdot; A_+, A_-)$ and  
$\xi\big(\, \cdot \, ; |\bsD_\bsA^{*}|^2, |\bsD_\bsA^{}|^2\big)$ (in the Lebesgue point sense).

\subsection{Computing the Witten index of the operator $\bsD_\bsA^{}$.} 
As a consequence of the principal trace formula, Theorem \ref{ptf}, and the Lifshitz--Krein trace formula,  the following equality holds, 
\begin{equation}\label{for_comput_of_WI}
z \tr_{L^2(\bbR;\cH)}\Big(\big(|\bsD_\bsA^{*}|^2 - z \, \bsI\big)^{-1} - 
\big(|\bsD_\bsA^{}|^2 - z \, \bsI\big)^{-1}\Big) = 
- \f{z}{2} \int_{\bbR} \frac{\xi(\nu; A_+, A_-) \, d\nu}{(\nu^2 - z)^{3/2}}, 
\end{equation}
 for all $z \in \bbC \backslash [0,\infty).$  
Recalling that the resolvent regularized Witten index $W_r(\bsD_\bsA^{})$ is the limit of the LHS as $z\rightarrow 0, \, z<0$, we see that this index can be computed by taking the limit of the RHS as $z\rightarrow 0, \, z<0$. To this end, we consider the operator $\mathbf{T}$, defined by setting
\begin{equation*} 
\mathbf{T}: \begin{cases} 
L^1\big(\bbR; (1+\nu^2)^{-3/2}d\nu\big) \to {\rm Hol} (\bbC \backslash [0,\infty))  \\[1mm] 
f \mapsto -z\int_{\mathbb{R}} d \nu \, (\nu^2-z)^{-3/2} f(\nu), \quad z \in \mathbb{C}\backslash [0,\infty),
\end{cases} 
\end{equation*} 
where ${\rm Hol} (\bbC \backslash [0,\infty))$ denotes the set of all holomorphic functions on $\bbC \backslash [0,\infty).$

\begin{lemma}\cite{CGPST15}\label{lem_4}
If $0$ is a left and a right Lebesgue point for $f\in L^1\big(\bbR; (1+\nu^2)^{-3/2}d\nu\big)$, then 
\begin{equation} 
\lim_{z\to0, \, z<0}(\mathbf{T}f)(z)=\Lf(0_+) + \Lf(0_-).
\end{equation} 
\end{lemma}

Thus, applying this lemma to the function $\xi(\, \cdot \,; A_+, A_-)$ on the right hand side of \eqref{for_comput_of_WI}, we arrive at the equality 
\begin{align*}
W_r(\bsD_\bsA^{})&=\lim_{z\to0, \, z<0}z \tr_{L^2(\bbR;\cH)} 
\Big(\big(|\bsD_\bsA^{*}|^2 - z \, \bsI\big)^{-1} - 
\big(|\bsD_\bsA^{}|^2 - z \, \bsI\big)^{-1}\Big)\\
&=\lim_{z\to0, \, z<0}- \f{z}{2} \int_{\bbR} \frac{\xi(\nu; A_+, A_-) \, d\nu}{(\nu^2 - z)^{3/2}}\\
&=\frac12\lim_{z\to0, \, z<0}(\mathbf{T}\xi(\, \cdot \,; A_+, A_-))(z)\stackrel{Lemma \ref{lem_4}}{=}[\Lxi(0_+) + \Lxi(0_-)]/2.
\end{align*}

Now, we turn to computing the semigroup regularized Witten index $W_s(\bsD_\bsA^{})$. 
To this end, we have established the following 
\begin{theorem}\cite{CGPST15}
If $0$ is a right Lebesgue point of $\xi\big(\, \cdot \,; |\bsD_\bsA^{*}|^2, |\bsD_\bsA^{}|^2\big)$, 
then
\begin{equation*}
\lim_{z \to \infty, z>0} \tr_{L^2(\bbR;\cH)} \Big(e^{-z  |\bsD_\bsA^{*}|^2} - e^{- z |\bsD_\bsA^{}|^2} \Big) 
= \Lxi\big(0_+; |\bsD_\bsA^{*}|^2, |\bsD_\bsA^{}|^2\big), 
\end{equation*}
uniformly with respect to $z$.
\end{theorem} 
Therefore, if $0$ is a right and a left Lebesgue point of $\xi(\, \cdot \,; A_+,A_-)$, then combining this theorem with equality \eqref{eq_in_Lebesgue_sense} we obtain 
$$W_s(\bsD_\bsA^{})=\Lxi\big(0_+; |\bsD_\bsA^{*}|^2, |\bsD_\bsA^{}|^2\big)\stackrel{Lemma \, \ref{lem_3}}{=}[\Lxi(0_+) + \Lxi(0_-)]/2.$$

\begin{theorem}\cite[Theorem 4.3]{CGPST15}\label{main}
 Assume Hypothesis \ref{hyp} and suppose 
that $0$ is a right and a left Lebesgue point of $\xi(\, \cdot \,;  A_+, A_-)$. Then $0$ is a right Lebesgue point of $\xi\big(\, \cdot \, ; |\bsD_\bsA^{*}|^2, |\bsD_\bsA^{}|^2\big)$ and
\begin{align*} 
\quad \Lxi\big(0_+; |\bsD_\bsA^{*}|^2, |\bsD_\bsA^{}|^2\big) 
= [\Lxi(0_+; A_+,A_-) + \Lxi(0_-; A_+, A_-)]/2.  
\end{align*}
and
\begin{align*} 
& W_r(\bsD_\bsA^{}) = [\Lxi(0_+; A_+,A_-) + \Lxi(0_-; A_+, A_-)]/2 = W_s(\bsD_\bsA^{}). 
\end{align*} 
\end{theorem} 

We emphasize that Theorem \ref{main} contains Theorem \ref{thm_Fredholm_case} as a particular case. Indeed, suppose that the operator $\bsD_\bsA^{}$ is Fredholm, that is, the asymptotes $A_\pm$ are boundedly invertible. In this case, $0$ is a right and a left Lebesgue point of $\xi(\, \cdot \,;  A_+, A_-)$ and  $[\Lxi(0_+; A_+,A_-) + \Lxi(0_-; A_+, A_-)]/2 =\xi(0; A_+,A_-)$. 

In the next subsection we discuss the case when $0$ may belong to the spectra of the operators $A_+$ and $A_-$. 
As we already noted the Witten index, in general, can be any prescribed 
real number. Next we demonstrate that this also applies to the special case of the Witten 
index of $\bsD_\bsA^{}$.

A simple concrete example is the following: Consider $A_{\pm} \in \cB(\cH)$ 
with $[A_+ - A_-] \in \cB_1(\cH)$, and introduce  the family 
\begin{equation*}
A(t) = A_-+{e^t}({e^t + 1})^{-1} [A_+ - A_-], \quad t \in \bbR,
\end{equation*} which satisfies Hypothesis \ref{hyp}. Moreover, since \emph{any} integrable function 
$\xi \in L^1(\bbR; dt)$ of compact support arises as the 
spectral shift function for a pair of bounded, self-adjoint operators 
$(A_+,A_-)$ in $\cH$ with $[A_+ - A_-] \in \cB_1(\cH)$, Theorem \ref{main} implies  that 
\begin{align*}
\begin{split}
W_r(\bsD_\bsA^{}) &= W_s(\bsD_\bsA^{})
= \Lxi(0; A_+,A_-)   \\
&= \text{ any prescribed real number.}   
\end{split}
\end{align*}

\subsection{The spectra of $A_\pm$ and Lebesgue points of $\xi(\, \cdot \,;  A_+,A_-)$.} 
We start with the simpler case where $A_{\pm}$ have discrete 
spectrum in an open neighbourhood of $0$. That is we assume, that for some $\nu>0$, the interval $(-\nu,\nu)$ contains only eigenvalues of $A_\pm$ of \emph{finite multiplicities}, which are \emph{isolated points} in $\sigma(A_\pm)$. The following remark easily follows from properties of SSF (see Subsection \ref{sec_properties_ssf}, property $(iii)$).

\begin{remark} 
If $A_{\pm}$ have discrete spectra in an open neighborhood 
of $0$, then the SSF $\xi(\, \cdot \,; A_+,A_-)$ has a right and left limit at any point of this open neighborhood 
and, in particular, any point in that open neighborhood is a right and a left Lebesgue point of 
$\xi(\, \cdot \,; A_+,A_-)$. \hfill $\diamond$
\end{remark}

On the contrary, in the presence of purely absolutely 
continuous spectrum of $A_{\pm}$ in a neighborhood of $0$, the situation is more complicated.

\begin{proposition}\cite[Proposition 4.6]{CGPST15}
There exist pairs of bounded self-adjoint operators $(A_+, A_-)$ in $\cH$ such that $(A_+ - A_-)$ 
is of rank-one, and $A_{\pm}$ both have purely absolutely continuous spectrum in a fixed 
neighborhood $(-\varepsilon_0, \varepsilon_0)$, for some $\varepsilon_0 > 0$, yet 
$\xi(\, \cdot \,; A_+,A_-)$ \emph{may or may not have} a right and/or a left Lebesgue point at $0$. 
\end{proposition}

\subsection{The Witten index of $\bsD_\bsA^{}$ in the Special Case $\dim(\cH) < \infty$.} 
We briefly treat the special finite-dimensional 
case, $\dim(\cH) < \infty$,  and explicitly compute the Witten index of 
$\bsD_\bsA^{}$ \emph{irrespectively} of whether or not $\bsD_\bsA^{}$ is a 
Fredholm operator in $L^2(\bbR; \cH)$. 

In this special case the Hypothesis \ref{hyp} acquires a considerably simpler form. We just suppose that  
\begin{equation} 
A_- \in \cB(\cH) \, \text{ is a self-adjoint matrix in $\cH$,}  \lb{eq26}
\end{equation}
and there exists a family of bounded self-adjoint matrices $\{A(t)\}_{t\in\bbR}$, locally absolutely continuous on $\bbR$, such that 
\begin{equation}  
\int_\bbR dt \, \big\|A'(t)\big\|_{\cB(\cH)} < \infty.   \lb{eq27}
\end{equation}

In the following we denote by 
$
\#_{>} (S) \, \text{ (respectively\ $\#_{<} (S)$)}
$
the number of strictly positive (respectively, strictly negative) eigenvalues of a 
self-adjoint matrix $S$ in $\cH$, counting multiplicity.   
Under these assumptions the formula for the Witten index of the operator $\bsD_\bsA^{}$ takes a particularly simple form. If should be pointed out that the result below yields (in a very special 
setting where $\dim(\cH)=1$) the result of Example \ref{exampleBGGSS}. 

\begin{theorem}\cite[Theorem 5.2]{CGPST15} \lb{t4.2} 
Assume Hypotheses \eqref{eq26} and \eqref{eq27}. Then the SSF $\xi(\,\cdot\,; A_+,A_-)$ has a
piecewise constant representative on $\bbR$, the right limit  
$\xi\big(0_+; |\bsD_\bsA^{*}|^2, |\bsD_\bsA^{}|^2\big)$ exists, 
and the SSF $\xi\big(\,\cdot\, ; |\bsD_\bsA^{*}|^2, |\bsD_\bsA^{}|^2\big)$ has a continuous 
representative on $(0,\infty)$. Moreover, the resolvent and semigroup regularized Witten 
indices $W_r(\bsD_\bsA^{})$ and $W_s(\bsD_\bsA^{})$ exist, and 
\begin{align*}
W_r(\bsD_\bsA^{}) = W_s(\bsD_\bsA^{}) &= \xi\big(0_+; |\bsD_\bsA^{*}|^2, |\bsD_\bsA^{}|^2\big)     \\
&= [\xi(0_+; A_+,A_-) + \xi(0_-; A_+,A_-)]/2     \\
&= \f{1}{2} [\#_{>} (A_+) - \#_{>} (A_-)] - 
\f{1}{2} [\#_{<} (A_+) - \#_{<} (A_-)].    
\end{align*}
In particular, in the finite-dimensional context, the Witten indices are  either  
integer, or half-integer-valued.
\end{theorem}

\section{Further extensions.}\label{dim_1} 
In this section we discuss an important example of operators $A_+$ and $A_-$, whose spectra are absolutely continuous and coincide with the whole real line and for which the results of previous sections are not applicable. 
The results of \cite{Pu08} and \cite{GLMST11}, \cite{CGPST15} describe  the Fredholm/Witten index theory for operators permitting essential spectra but the relatively trace class  assumption rules out standard partial differential operators such as Dirac type operators.
Thus, in order to incorporate this important class of examples, we need a more general framework.

To illustrate this fact, we consider the following example.
Let $A_-$ and $\{A(t)\}_{t \in \bbR}$ be given by 
$$A_-=\frac{d}{idx}, \quad A(t)=A_-+\theta(t)M_f,\quad \dom(A_-)=\dom(A(t))=W^{1,2}(\bbR),
\quad t \in \bbR,$$
that is, we consider  the differentiation operator on $L^2(\bbR;dx)$ and its perturbation by multiplication operator $M_f$ defined by a function  $f\in L^\infty(\bbR;dx)$.
Here   $\theta$ is a function satisfying 
\begin{align*}
& 0 \leq \theta \in L^{\infty}(\bbR; dt),\; \;
\theta' \in L^{\infty}(\bbR; dt) \cap L^1(\bbR; dt),\\ 
& \lim_{t\rightarrow-\infty}\theta(t)=0, \quad  \lim_{t\rightarrow+\infty}\theta(t)=1.
\end{align*}

Then the asymptotes $A_\pm$ of the family $\{A(t)\}_{t\in\bbR}$ as $t\rightarrow\pm\infty$ are given by $A_-$ and 
$$A_+=A_-+M_f.$$
In other words, we have a one-dimensional Dirac operator and its perturbation by a bounded function. The well-known Cwikel's estimates (see, e.g., \cite[Ch.~4]{Si05}) guarantee that for $f$  decaying sufficiently rapidly at $\pm\infty$, the operator 
$(A_+-A_-)(A_-^2 + 1)^{-s/2}$ is trace class for $s>1$, but for no lesser 
value of $s$. Thus, even in one dimension, the relative trace class assumption is violated for the example above.

However, although the one-dimensional differential operator $A_-$ and its perturbations 
do not satisfy the relative trace class assumption, we still can compute the Witten index 
$W_r(\bsD_\bsA^{})$. 
For this one-dimensional setting, under the identification of the Hilbert spaces 
$L^2(\bbR; dt; L^2(\bbR; dx))$ and $L^2(\bbR^2; dtdx)$,
the operator $\bsD_\bsA^{}$, defined by \eqref{def_D_A}, is given by
$$\bsD_\bsA^{}=\frac{d}{dt}+\bsA,$$
with $\bsA=\frac{d}{idx}+M_\theta M_f.$ That is, in this setting we work with the operator 
$$\bsD_\bsA^{}=\frac{d}{dt}+\frac{d}{idx}+M_\theta M_f.$$

Since the operator $\frac{d}{idx}$ has absolutely continuous spectrum, coinciding with the whole real line, the operator $\bsD_\bsA^{}$ possesses the following properties: \\[1mm] 
$(i)$ Since $0\in\sigma(A_-)=\bbR$, by Theorem \ref{thn_Fredholm} we have that the operator $\bsD_\bsA^{}$  is \emph{not Fredholm.} \\[1mm] 
$(ii)$ The essential spectrum of the operator $\bsD_\bsA^{}$ is the \emph{whole complex plane $\bbC$} (see Corollary \ref{spec_ess}).

It is interesting (and somewhat surprising) that for this particular example under some assumptions on the perturbation $f$ (see Theorem \ref{main2} below) we still have the inclusions (cf. Theorem \ref{ptf})
\begin{align*} 
[g_z(A_+)-g_z(A_-)]& \in\cB_1(L^2(\bbR)),\\
\Big(\big(|\bsD_\bsA^{*}|^2 - z \, \bsI\big)^{-1} - 
\big(|\bsD_\bsA^{}|^2 - z \, \bsI\big)^{-1}\Big) 
&\in \cB_1\big(L^2\big(\bbR^2)\big), \quad z\in\bbC \backslash [0,\infty),
\end{align*} where 
 $g_z(x)=x(x^2-z)^{-1/2}$, $x \in \bbR$. Moreover, using an approximation technique, we can prove the principal trace formula as in Theorem \ref{ptf}
 \begin{align*} &\tr_{L^2(\bbR^2)}
\Big(\big(|\bsD_\bsA^{*}|^2 - z \, \bsI\big)^{-1} - 
\big(|\bsD_\bsA^{}|^2 - z \, \bsI\big)^{-1}\Big) = \frac{1}{2z}\tr_{L^2(\bbR)} (g_z(A_+)-g_z(A_-)),  
\end{align*}
for all $z\in\bbC \backslash [0,\infty)$.

The main application of this principal trace formula is an extension of Pushnitski's formula. 
 Furthermore, employing some classical harmonic analysis we are able to compute the \emph{actual value} of (a representative of) the spectral shift function for the pair $A_+,A_-$. 

\begin{theorem}\label{main2}
Let $
f \in W^{1,1}(\bbR; dx)\cap C_b(\bbR; dx)$ and $ f'\in L^\infty(\bbR;dx).$ Then for a.e.~$\lambda>0$ 
and a.e.~$\nu \in \bbR$, 
$$
\xi\big(\lambda; |\bsD_\bsA^{*}|^2, |\bsD_\bsA^{}|^2\big)=\xi(\nu; A_+,A_-) 
= \frac{1}{2\pi}\int_\bbR f(x)\,dx. 
$$ 
\end{theorem}

The fact that the SSF $\xi(\, \cdot \,; A_+,A_-)$ is a constant immediately implies that $0$ is Lebesgue point of the function $\xi(\, \cdot \,; A_+,A_-)$. 

\begin{theorem}\cite{CGLPSZ15}
The Witten indices $W_r(\bsD_\bsA^{})$ and $W_s(\bsD_\bsA^{})$ of the operator  $\bsD_\bsA^{}$ exist and equal 
\begin{equation*}
W_r(\bsD_\bsA^{}) = W_s(\bsD_\bsA^{}) = \xi\big(0_+;|\bsD_\bsA^{*}|^2, |\bsD_\bsA^{}|^2\big) = 
\xi(0; A_+, A_-) = \f{1}{2 \pi} \int_{\bbR} f(x) dx. 
\end{equation*} 
\end{theorem}

\begin{remark}
We note that the equality $\xi(\, \cdot \,; A_+, A_-) = \f{1}{2 \pi} \int_{\bbR}f(x) dx$ may  also be proved via scattering theory and modified Fredholm determinants of 2nd order 
(cf.\ \cite{CGLPSZ15}). \hfill $\diamond$
\end{remark}

The results above can be also generalized to the following setting.
Assume that $A_-$ is an (unbounded) self-adjoint operator in a complex separable Hilbert space $\cH$ and assume that the family of bounded operators $\{B(t)\}_{t\in\bbR}$ is a $2$-relative trace class perturbation, that is, $B'(t)(|A_-| + 1)^{-2}\in\cB_1(\cH)$,  $t\in\bbR$, and
\begin{equation*}  
\int_\bbR dt \, \big\|B'(t) (|A_-| + 1)^{-2}\big\|_{\cB_1(\cH)} < \infty.
\end{equation*}Imposing some minor additional conditions on the family $\{B(t)\}_{t\in\bbR}$ one can prove the following result:

\begin{theorem}
Suppose that $0$ is a right 
and a left Lebesgue point of $\xi(\,\cdot\,\, ; A_+, A_-)$, then $0$ is 
also a right Lebesgue point of $\xi\big(\,\cdot\,\, ;|\bsD_\bsA^{*}|^2, |\bsD_\bsA^{}|^2\big)$ 
and $W_r(\bsD_\bsA^{})$ exists and equals 
\begin{equation*}
W_r(\bsD_\bsA^{}) = \Lxi\big(0_+;|\bsD_\bsA^{*}|^2, |\bsD_\bsA^{}|^2\big) 
= [\Lxi(0_+; A_+,A_-) + \Lxi(0_-; A_+, A_-)]/2.    
\end{equation*}
\end{theorem}

\medskip 
 
\noindent 
{\bf Acknowledgments.} A.C., G.L., and F.S. gratefully acknowledge financial support from the Australian Research Council. 

F.S. sincerely thanks the organizers of the conference, {\it Spectral Theory and Mathematical Physics,} held in Santiago, Chile, in November of 2014, for the hospitality extended to him and for 
the opportunity to deliver a mini-course on the subject matter treated in this survey. The present work is a substantially revised and extended version of that course.



\begin{thebibliography}{999}

\bibitem{AC79} Y.\ Aharonov and A.\ Casher, {\it A. Ground state of a spin-$1/2$ 
charged particle in a two-dimensional magnetic field}, Phys. Rev. A (3) {\bf 19}, 
2461--2462 (1979).

\bibitem{An89} N.\ Anghel, {\it Remark on Callias' index theorem}, Rep. Math. Phys. 
{\bf 28}, 1--6 (1989).

\bibitem{An90} N.\ Anghel, {\it $L^2$-index formulae for perturbed Dirac operators},  
Comm. Math. Phys. {\bf 128}, 77--97 (1990).

\bibitem{An90a} N.\ Anghel, {\it The two-dimensional magnetic field problem 
revisited}, J. Math. Phys. {\bf 31}, 2091--2093 (1990).

\bibitem{An93} N.\ Anghel, {\it On the index of Callias-type operators}, Geom. 
Funct. Anal. {\bf 3}, 431--438 (1993).

\bibitem{An94} N.\ Anghel, {\it Index theory for short-ranged fields in higher 
dimensions}, J. Funct. Anal. {\bf 119}, 19--36 (1994).

\bibitem{APS76}
M.\ F.\ Atiyah, V.\ K.\ Patodi, and I.\ M.\ Singer, {\it Spectral asymmetry and 
Riemannian geometry. III}, Math. Proc. Cambridge Philos. Soc. {\bf 79} 71--99 (1976).

\bibitem{ACDS09}  N.\ A.\ Azamov, A.\ L.\ Carey, P.\ G.\ Dodds, and F.\ A.\ Sukochev, 
{\it Operator integrals, spectral shift, and spectral flow}, Canad. J. Math. {\bf 61}, 
241--263 (2009).

\bibitem{ACS07} N.\ A.\ Azamov, A.\ L.\ Carey, and F.\ A.\ Sukochev, {\it The spectral shift
 function and spectral flow}, Comm. Math. Phys. {\bf 276}, 51--91 (2007).
 
\bibitem{ADS06} N.\ A.\ Azamov, P.\ G.\ Dodds, F.\ A.\ Sukochev, {\it The Krein 
spectral shift function in semifinite von Neumann algebras.} Integral Eqs. Operator 
Theory {\bf 55}, 347--362 (2006).

\bibitem{BS75} M.\ Sh.\ Birman and M.\ Z.\ Solomyak, {\it Remarks on the 
spectral shift function}, J. Sov. Math. {\bf 3}, 408--419 (1975).  

\bibitem{BS78} R.\ Bott and R.\ Seeley, {\it Some remarks on the paper of Callias}, 
Comm. Math. Phys. {\bf 62}, 235--245 (1978).

\bibitem{BC15} M.\ Braverman and S.\ Cecchini, {\it Spectral theory of von Neumann algebra valued differential operators over non-compact manifolds}, arXiv{1503.02998}.

\bibitem{Br86} L.\ G.\ Brown, {\it Lidski's theorem in the type II case.} Geometric Methods 
in Operator Algebras (Kyoto, 1983), H.\ Araki and E.\ G.\ Effros (eds.), Pitman Res. Notes 
Math. Ser., {\bf 123}, Longman Sci. Tech., Harlow, 1986, pp.\ 1--35.

\bibitem{BGGSS87} D.\ Boll{\' e}, F.\ Gesztesy, H.\ Grosse, W.\ Schweiger, and 
B.\ Simon, {\it Witten index, axial anomaly, and Krein's spectral 
shift function in supersymmetric quantum mechanics}, J. Math. Phys.
{\bf 28}, 1512--1525 (1987).

\bibitem{BMS88} N.\ V.\ Borisov, W.\ M\"uller, and R.\ Schrader, {\it Relative 
index theorems and supersymmetric scattering theory}, Commun. Math. 
Phys. {\bf 114}, 475--513 (1988).  

\bibitem{BB84} D.\ Boyanovsky and R.\ Blankenbecler, {\it Fractional
indices in supersymmetric theories}, Phys. Rev. {\bf D 30}, 1821--1824
(1984). 

\bibitem{Bu92} U.\ Bunke, {\it Relative index theory}, J. Funct. Anal. {\bf 105}, 
63--76 (1992).

\bibitem{Bu95} U.\ Bunke, {\it A K-theoretic relative index theorem and 
Callias-type Dirac operators}, Math. Ann. {\bf 303}, 241--279 (1995). 

\bibitem{Ca78} C.\ Callias, {\it Axial anomalies and index theorems on 
open spaces},  Commun. Math. Phys. {\bf 62}, 213--234 (1978).

\bibitem{CGLPSZ15} A.\ Carey, F.\ Gesztesy, G.\ Levitina, D.\ Potapov, F.\ Sukochev, 
and D.\ Zanin, {\it On index theory for non-Fredholm operators: a
1 + 1-dimensional example}, to appear. 

\bibitem{CK15} A.\ Carey and J.\ Kaad, {\it Topological invariance of the 
homological index}, J. reine angew. Math., to appear,  \arxiv{1402.0475}.

\bibitem{CGPST14} A.\ Carey, F.\ Gesztesy, D.\ Potapov, F.\ Sukochev, and
Y.\ Tomilov, {\it A Jost--Pais-type reduction of Fredholm determinants for 
semi-separable operators in infinite dimensions and some applications}, 
Integral Eqs. Operator Theory {\bf 79}, 389--447 (2014). 

\bibitem{CGPST15} A.\ Carey, F.\ Gesztesy, D.\ Potapov, F.\ Sukochev, and
Y.\ Tomilov, {\it On the Witten Index in terms of spectral shift functions}, J. Analyse Math., 
to appear, \arxiv{1404.0740}.  

\bibitem{CP86} R.\ W.\ Carey and J.\ D.\ Pincus, {\it Index theory for operator ranges and 
geometric measure theory}, in {\it Geometric Measure Theory and the Calculus of Variations 
(Arcata, Ca, 1984)}, W.\ K.\ Allard and F.\ J.\ Almgren (eds.), Proc. Symposia in Pure Math., 
Vol.\ 44, Amer. Math. Soc., Providence, RI, 1986, pp.\ 149--161.

\bibitem{EE89} D.\ E.\ Edmunds and W.\ D.\ Evans, {\it Spectral Theory and
Differential Operators}, Clarendon Press, Oxford, 1989.

\bibitem{EGH80} T.\ Eguchi, P.\ B.\ Gilkey, and A.\ J.\ Hanson, {\it Gravitation, gauge 
theories and differential geometry}, Phys. Rep. {\bf 66}, No.\ 6, 213--393 (1980). 

\bibitem{FOW87} P.\ Forgacs, L.\ O'Raifeartaigh, and A.\ Wipf, {\it Scattering theory, $U(1)$ 
anomaly and index theorems for compact and non-compact manifolds}, Nuclear Phys. 
{\bf B 293}, 559--592 (1987). 

\bibitem{Ge86} F.\ Gesztesy, {\it Scattering theory for one-dimensional 
systems with nontrivial spatial asymptotics},  in {\it Schr\"odinger operators, 
Aarhus 1985}, Lecture Notes in Math., Vol.\ 1218, Springer, Berlin, 1986, pp.\ 93--122.

\bibitem{GLMST11} F.\ Gesztesy, Yu.\ Latushkin, K.\ Makarov, F.\ Sukochev and Yu.\ Tomilov, {\it
The index formula and the spectral shift function for relatively trace class perturbations.} 
Adv. Math. {\bf 227}, 319--420 (2011).

\bibitem{GPS08} F.\ Gesztesy, A.\ Pushnitski, and B.\ Simon, {\it On the Koplienko spectral 
shift function. I. Basics}, J. Math. Phys., Anal., Geometry, {\bf 4}, No.\ 1, 63--107 (2008). 

\bibitem{GS88} F.\ Gesztesy and B.\ Simon, {\it Topological invariance of the Witten index}, 
J. Funct. Anal. {\bf 79}, 91--102 (1988). 

\bibitem{GK69} I.\ Gohberg and M.\ G.\ Krein, {\it Introduction to the Theory of
Linear Nonselfadjoint Operators}, Translations of Mathematical Monographs,
Vol.\ 18, Amer. Math. Soc., Providence, RI, 1969.

\bibitem{HS65} E.\ Hewitt and K.\ Stromberg, {\it Real and Abstract
Analysis. A Modern Treatment of the Theory of Functions of a Real Variable},
Springer, New York, 1965.

\bibitem{Hi83} M.\ Hirayama, {\it Supersymmetric quantum mechanics
and index theorem}, Progr. Theoret. Phys. {\bf 70}, 1444--1453 (1983). 

\bibitem{Hi86} M.\ Hirayama, {\it Topological invariants for Dirac operators on open 
spaces}, Phys. Rev. {\bf D 33}, 1079--1087 (1986).  

\bibitem{IM84} C.\ Imbimbo and S.\ Mukhi, {\it Topological invariance in supersymmetric 
theories with a continuous spectrum}, Nuclear Phys. {\bf B 242}, 81--92 (1984).  

\bibitem{Ka80} T.\ Kato, {\it Perturbation Theory for Linear Operators}, corr.\ printing 
of the 2nd ed., Springer, Berlin, 1980.

\bibitem{Ko84}  L.\ S.\ Koplienko, {\it Trace formula for nontrace-class perturbations}, 
Sib. Math. J. {\bf 25}, 735--743 (1984).

\bibitem{Ko86} L.\ S.\ Koplienko, {\it Local conditions for the existence of the 
spectral shift function}, J. Sov. Math. {\bf 34}, 2080--2090 (1986). 

\bibitem{Ko11} C.\ Kottke, {\it An index theorem of Callias type for pseudodifferential operators}, 
J. K-Theory {\bf 8}, no.\ 3, 387--417 (2011).

\bibitem{Ko15} C.\ Kottke, {\it A Callias-type index theorem with degenerate potentials}, 
Commun. PDE {\bf 40}, 219--264 (2015).

\bibitem{Kr53}  M.\ G.\ Krein, {\it On the trace formula in perturbation theory} (Russian), 
Mat. Sbornik N.S. {\bf 33} (75), 597--626 (1953). 

\bibitem{Kr62}  M.\ G.\ Krein, {\it Perturbation determinants and a formula for the traces of 
unitaryand self-adjoint operators}, Soviet Math.\ Dokl.\ {\bf 3}, 707--710 (1962).

\bibitem{KY81} M.\ G.\ Krein and V.\ A.\ Yavryan, {\it Spectral shift functions that arise in 
perturbations of a positive operator} (Russian), J. Operator Theory {\bf 6}, 155--191 
(1981).

\bibitem{Mu87} W.\ M\"uller, {\it Manifolds with Cusps of Rank One. Spectral Theory 
and $L^2$-Index Theorem}, Lecture Notes in Math., Vol.\ 1244, Springer, 
Berlin, 1987. 

\bibitem{Mu88} W.\ M\"uller, {\it $L^2$-index and resonances}, in {\it Geometry 
and Analysis on Manifolds}, T.\ Sunada (ed.), Lecture Notes in Math., Vol.\ 1339, Springer, Berlin, 1988, pp.\  203--221. 

\bibitem{Mu98} W.\ M\"uller, {\it Relative zeta functions, relative determinants and
scattering theory}, Commun. Math. Phys. {\bf 192}, 309--347 (1998).

\bibitem{NS86} A.\ J.\ Niemi and G.\ W.\ Semenoff, {\it Index theorems on open infinite manifolds}, 
Nuclear Phys. B  {\bf 269}, 131--169 (1986).

\bibitem{NS86a} A.\ J.\ Niemi and G.\ W.\ Semenoff, {\it Fermion number fractionization 
in quantum field theory}, Phys. Rep. {\bf 135}, No.\ 3, 99--193 (1986). 

\bibitem{Ph96} J.\ Phillips, {\it Self-adjoint Fredholm operators and spectral flow}, 
Canad. Math. Bull. {\bf 39}, 460--467 (1996). 

\bibitem{PSS13} D.\ Potapov, A.\ Skripka, and F.\ Sukochev, {\it Spectral shift function 
of higher order}, Invent. Math. {\bf 193}, 501--538 (2013).

\bibitem{PSZ} D.\ Potapov, F.\ Sukochev, D.\  Zanin,  {\it  Krein's trace theorem revisited}, 
J. Spectr. Theory {\bf 4}, 415--430 (2014).

\bibitem{Pu08} A.\ Pushnitski, {\it The spectral flow, the Fredholm index, and the spectral 
shift function}, in {\it Spectral Theory of Differential Operators: M.\ Sh.\ 
Birman 80th Anniversary Collection}, T.\ Suslina and D.\ Yafaev (eds.),
AMS Translations, Ser.\  2, Advances in the Mathematical Sciences, Vol.\  225,
Amer. Math. Soc., Providence, RI, 2008, pp.\ 141--155.

\bibitem{RS95} J.\ Robbin and D.\ Salamon, {\it  The spectral flow and the Maslov
index},  Bull. London Math. Soc. {\bf 27}, 1--33 (1995).

\bibitem{Si05} B.\ Simon, {\it Trace Ideals and Their Applications}, 2nd ed.,
Mathematical Surveys and Monographs, Vol.\ 120, Amer. Math. Soc.,
Providence, RI, 2005.

\bibitem{Si94} K.\ Sinha K. and A.\ Mohapatra, {Spectral shift function and trace
formula}, Proc. Indian Acad. Sci. {\bf 104}, 819--853 (1994).

\bibitem{So93} A.\ V.\ Sobolev, {\it Efficient bounds for the spectral shift function}, 
Ann. Inst. H. Poincar\'e {\bf A 58}, 55--83 (1993). 

\bibitem{Th92} B.\ Thaller, {\it The Dirac Equation}, Texts and Monographs in 
Physics, Springer, Berlin, 1992. 

\bibitem{Vo87} D.\ Voiculescu, {\it On a trace formula of M. G. Krein,} in 
{\it Operators in Indefinite Metric Spaces, Scattering Theory and Other Topics, 
(Bucharest, 1985)}, H.\ Helson and G.\ Arsene (eds.), Oper. Theory Adv. Appl., 
Vol.\ 24, Birkhauser, Basel, 1987, pp.\ 329--332.

\bibitem{Wi82} E.\ Witten, {\it Constraints on supersymmetry breaking}, Nuclear Phys.
{B 202}, 253--316 (1982).

\bibitem{Ya92} D.\ R.\ Yafaev, {\it Mathematical Scattering Theory.
General Theory}, Transl. Math. Monographs, Vol.\ 105, Amer. Math. Soc., Providence, RI, 
1992.

\bibitem{Ya05} D.\ R.\ Yafaev, {\it A trace formula for the Dirac operator}, 
Bull. London Math. Soc. {\bf 37}, 908--918 (2005).

\end{thebibliography}
\end{document}